\documentclass[12pt]{amsart}

\usepackage{amssymb,amsmath,amsthm}
\usepackage{graphics,graphicx}
\usepackage{tikz-cd}
\usetikzlibrary{calc}
\usepackage{dynkin-diagrams}
\usepackage{dsfont}
\usepackage{xcolor}
\usepackage{enumitem}
\usepackage{booktabs}

\usepackage[a4paper]{geometry}
\usepackage[colorlinks=true,allcolors=blue]{hyperref} %

\newtheorem{thm}{Theorem}[section]
\newtheorem{cor}[thm]{Corollary}
\newtheorem{lem}[thm]{Lemma}
\newtheorem{prop}[thm]{Proposition}

\theoremstyle{definition}
\newtheorem{defn}[thm]{Definition}
\newtheorem{rem}[thm]{Remark}
\newtheorem{exa}[thm]{Example}

\numberwithin{equation}{section} \numberwithin{thm}{section}

\begin{document}

\title[Spherical secondary polytopes]{Spherical secondary polytopes and equivariant Chow polytopes}

\author{Thibaut Delcroix}
\address{Thibaut Delcroix, IMAG, Univ Montpellier, CNRS, Montpellier, France}
\email{thibaut.delcroix@umontpellier.fr}
\urladdr{http://delcroix.perso.math.cnrs.fr/}

\author{King Leung Lee}

\author{Naoto Yotsutani}
\address{Naoto Yotsutani, Department of Mathematics, Faculty of Science, Shizuoka University, $836$ Ohya, Suruga-ku, Shizuoka-shi, Shizuoka, $422$-$8529$, Japan}
\email{yotsutani.naoto@shizuoka.ac.jp}

\address{IMAG, Univ Montpellier, CNRS, Montpellier, France}
\email{naoto.yotsutani@umontpellier.fr}

\date{\today}

\begin{abstract}
We introduce the spherical secondary polytope of a polarized spherical variety, defined as a weighted Minkowski integral, and identify it with its equivariant Chow polytope, generalizing the classical result for toric varieties.
We deduce convex criteria for Chow semistability and polystability.
For rank-one spherical varieties, we obtain an effective criterion for asymptotic Chow polystability in the K-polystable regime.
We also describe equivariant degenerations of spherical varieties in a fixed projective embedding in terms of satellites.
\end{abstract}

\maketitle


\section{Introduction}

Geometric invariant theory applies straightforwardly to projective hypersurfaces by encoding hypersurfaces by their equations. 
For a general projective subvariety, the Chow form plays the role of a defining equation. Applying geometric invariant theory to its projective class gives Chow stability~\cite{Mumford_1977}.
The same construction applies to effective cycles.
For smooth polarized varieties, Chow polystability is characterized by the existence of a balanced metric, providing a K\"ahler-geometric interpretation of Chow stability.
This criterion goes back to Zhang and Luo~\cite{Zhang_1996,Luo_1998}.

In the spirit of geometric quantization, Donaldson proved that existence of balanced metrics is closely related with existence of constant scalar curvature K\"ahler metrics: on a smooth polarized manifold \((X,L)\) with discrete automorphism group which admits a cscK metric in \(c_1(L)\), balanced metrics exist for the Kodaira embeddings associated with \(L^{\otimes k}\) for \(k\) large enough, and these metrics converge to the cscK metric~\cite{Donaldson_2001}. 
The relation between Chow forms and K\"ahler energy functionals was studied by Phong and Sturm~\cite{Phong_Sturm_2003} and Paul~\cite{Paul_2004}. 
For the algebro-geometric counterpart of the existence of cscK metrics, namely K-stability, this reflects a relationship between (asymptotic) Chow stability and K-stability, which was first explored by Ross and Thomas~\cite{Ross-Thomas_2007}. 

It should be highlighted that these relations are in fact quite subtle. 
Ono, Sano and the third named author~\cite{Ono-Sano-Yotsutani_2012} exhibited examples of smooth Fano manifolds with K\"ahler-Einstein metrics, hence K-polystable, but which are not Chow semistable for the Kodaira embedding associated with \(L^{\otimes k}\) for infinitely many values of \(k\). 
The phenomenon exhibited by this example was later clarified by Mabuchi: he proved that a smooth polarized manifold \((X,L)\)
admitting a cscK metric in \(c_1(L)\) is asymptotically Chow
polystable if his obstruction to asymptotic Chow semistability
vanishes~\cite{Mabuchi_2004,Mabuchi_2005}. 

For singular varieties, K-polystability and vanishing of Mabuchi's obstruction can coexist with failure of asymptotic Chow semistability. 
The second author gives such examples among Gorenstein toric Fano varieties in every dimension \(n\geq6\)~\cite[Example~3.9]{Lee_2025}.

The complete understanding of Chow stability for toric varieties was essential in these developments. 
Although the picture was clarified only gradually after Ono's work~\cite{Ono_2011,Ono_2013}, a combinatorial necessary-and-sufficient criterion of Chow polystability for toric varieties follows from the foundational work of Kapranov, Sturmfels and Zelevinsky~\cite{Kapranov_Sturmfels_Zelevinsky_1992}, also recorded in the famous book~\cite{GKZ}. 
They showed that the Chow polytope of a toric variety coincides with its secondary polytope.

Spherical varieties extend the toric setting to actions of connected reductive groups. 
A normal \(G\)-variety is spherical if a Borel subgroup of \(G\) has an open orbit. 
For a projective spherical variety, the space of global sections of each power of a \(G\)-linearized ample line bundle is a multiplicity-free \(G\)-representation. 
In the toric case its irreducible summands are character lines, and in general they may have higher dimension.

In this article, we define a secondary polytope for projective spherical varieties using their combinatorial data, namely the moment polytope and the valuation cone. 
This spherical secondary polytope fully encodes Chow polystability of the variety and may be formulated in a way that is easily comparable with the K-polystability criterion of the first author~\cite{Delcroix_2023_KSSV2}. 

Let \(X\) be a projective spherical \(G\)-variety and \(L\) a very ample \(G\)-linearized line bundle. 
Write \(\Delta\) for its moment polytope, \(\mathcal V\) for the valuation cone of \(X\),
and \(d\mu_{\mathrm{DH}}\) for the Duistermaat--Heckman measure on \(\Delta\). 
Write \(H^0(X,L)=\bigoplus_{\lambda\in\Lambda_1}W_{1,\lambda}\), where \(W_{1,\lambda}\) is the irreducible summand of highest weight \(\lambda\). 
We define the \emph{Chow weight functional} by
\[
\mathcal C_1(f)
=
\sum_{\lambda\in\Lambda_1}
       f(\lambda)\frac{\dim W_{1,\lambda}}
     {h^0(X,L)}
-
\int_\Delta f\,\frac{d\mu_{\mathrm{DH}}}
     {\int_\Delta d\mu_{\mathrm{DH}}}.
\]

\begin{thm}[Chow stability criterion]
\label{main_thm_Chow_stability}
For the Kodaira embedding defined by \(L\), the variety \(X\) is Chow semistable if and only if
\[
\mathcal C_1(f)\geq 0
\]
for every rational convex piecewise-affine function \(f\) on \(\Delta\) whose affine pieces have slopes in \(-\mathcal V\).
It is Chow polystable if and only if these inequalities hold and \(\mathcal C_1(f)=0\) if and only if \(f\) is affine with slope in the linear part \(\operatorname{Lin}(\mathcal{V}) = \mathcal V\cap(-\mathcal V)\).
\end{thm}

The spherical secondary polytope for \((X,L)\) is defined as follows. 
Let \(\Sigma\) be the spherical roots of \(X\). 
For \(x\in\Delta\), set
\[
K_x=\left\{r\in\mathbb R_{\geq0}^{\Lambda_1}\ \middle|\
\sum_{\lambda\in\Lambda_1}r_\lambda=1,\quad
\sum_{\lambda\in\Lambda_1}r_\lambda\lambda\in x+\operatorname{Cone}(\Sigma) \right\}.
\]

\begin{defn}
The \emph{spherical secondary polytope} of \((X,L)\) is the Minkowski integral 
\[
\Sigma_{\mathrm{sec}}^G(X,L)
=(\dim X+1)!\int_\Delta K_x\,d\mu_{\mathrm{DH}}(x).
\]
\end{defn}

This Minkowski integral is a weighted version of the fiber-polytope construction of Billera and Sturmfels~\cite{Billera_Sturmfels_1992}, with respect to the Duistermaat--Heckman measure \(\mu_{\mathrm{DH}}\). 
This interpretation also suggests a practical way to compute it, see Remark~\ref{rem:BS}. 

The relation between the spherical secondary polytope and Chow stability is obtained via the notion of equivariant Chow polytope for \(G\)-sub-varieties in the projectivization of a multiplicity-free \(G\)-representation, that we introduce here. 

Let \(V=\bigoplus_{i\in I}V_i\) be a multiplicity-free \(G\)-representation and let \(Z\) be a \(G\)-stable effective cycle in \(\mathbb P(V)\). 
The centralizer \(\widetilde T_0=Z_{\operatorname{GL}(V)}(G)\) is the torus acting by independent scalars on the summands \(V_i\).
We define the \emph{equivariant Chow polytope} by 
\[
\operatorname{Chow}_{\mathrm{GL}}^G(Z)
=
\operatorname{conv}\bigl\{
  \chi\in X^*(\widetilde T_0)\mid R_{Z,\chi}\neq 0
\bigr\}
\subset X^*(\widetilde T_0)\otimes_{\mathbb Z}\mathbb R,
\]
where \(R_{Z,\chi}\) is the \(\chi\)-weight component of the Chow form of \(Z\).
The polytope \(\operatorname{Chow}_{\mathrm{SL}}^G(Z)\) is its image under restriction of characters to \(T_0=Z_{\operatorname{SL}(V)}(G)^\circ\).

\begin{thm}[Equivariant Chow polytope formula]
\label{main_thm_Chow_polytope}
For the Kodaira embedding \(X\subset\mathbb P(H^0(X,L)^*)\), under the natural identification of the real character space of \(\widetilde T_0\) with \(\mathbb R^{\Lambda_1}\),
\[
\operatorname{Chow}_{\mathrm{GL}}^G(X)
=
\Sigma_{\mathrm{sec}}^G(X,L).
\]
\end{thm}

The full statement in the article, Theorem~\ref{thm:spherical-secondary-Chow}, also treats a base-point-free \(G\)-invariant linear system \(E\subset H^0(X,L)\) defining a finite morphism \(\pi:X\to\mathbb P(E^*)\). In this setting the same formula identifies the equivariant Chow polytope of the push-forward cycle \(\pi_*[X]\) with the spherical secondary polytope \(\Sigma_{\mathrm{sec}}^G(X,L;E)\) associated with \(E\).

This finite-morphism setting fits into the framework of Alexeev and Brion~\cite{Alexeev_Brion_2006}, who constructed a moduli space of stable spherical varieties over a fixed projective space.
Their cycle map sends a pair \((X,\pi)\) to the push-forward cycle \(\pi_*[X]\). 
In the toric case, the normalization of the image of the main part of their moduli space under this map is the toric variety associated with the classical secondary polytope.
For a fixed spherical pair \((X,\pi)\), our formula computes the weight polytope of the Chow form of \(\pi_*[X]\) under the torus that scales the irreducible summands of the ambient representation. 
We also describe the limiting cycles under one-parameter subgroups of this torus in terms of the spherical envelope, introduced here, and the satellites of the open orbit~\cite{Batyrev_Moreau_2020}. 

Although the spherical secondary polytope is fully encoded by the combinatorial data encoding spherical varieties, the flavor is somewhat different from the toric setting. 
The roles of the Weyl dimension formula and the valuation cone make this description more remote from triangulation-based characterization of the secondary polytope in the toric case.
In the horospherical case, vertices of the spherical secondary polytope are still in bijection with regular triangulations of the moment polytope, but the coordinates involve the Duistermaat-Heckman polynomial volumes. 
In the non-horospherical case, that is, when the valuation cone is not trivial, triangulations are no longer the meaningful decompositions in general, and maximal admissible decompositions of the polytope may further correspond to faces of positive dimension instead of vertices. 

We illustrate the construction with several examples, including quadrics equipped with several spherical actions. First, consider the smooth quadric fourfold
\(
X=\{q_U(u)+q_W(w)=0\}
\subset\mathbb P(U\oplus W),
\)
where \(U\) and \(W\) are three-dimensional quadratic spaces, equipped with the action of \(G=\mathrm{SO}(U)\times\mathrm{SO}(W)\). For \(L=\mathcal O_X(1)\), in natural coordinates we obtain
\(
\operatorname{Chow}_{\mathrm{GL}}^G(X)
=\operatorname{conv}\{(4,6),(6,4)
\}.
\)
The two vertices correspond to the same trivial decomposition of the moment polytope, which is a two dimensional simplex with only two integral points, but to distinct degenerations of \(X\). 
Thus, unlike in the toric case, the decomposition of the moment polytope alone does not determine the corresponding Chow vertex. 
See Examples~\ref{exa:two-block-quadric}, \ref{exa:quadric-same-decomposition} and~\ref{exa:GL-quadric-boundary} for further illustrations with quadrics.

Our third main result is a finite test for asymptotic Chow stability in spherical rank one, under a K-polystability assumption. 
In spherical rank one, we first obtain a finite test at a fixed level \(k\): it suffices to test an explicit finite list of convex piecewise affine functions, each with a single breakpoint at a normalized weight, together with an affine vanishing condition in the horospherical case. 
The corresponding
spherical envelopes have at most two domains of affinity, so each limit cycle has at most two irreducible components.
This is reminiscent of the work of Apostolov, Calderbank, Gauduchon and T{\o}nnesen-Friedman~\cite{ACGT_2008}, who use test configurations with two-component central fibers to test relative K-polystability of admissible projective bundles, as well as the celebrated work of Donaldson~\cite{Donaldson_2002} on K-stability of toric surfaces.
For the broader question of making stability criteria effective, see also the survey~\cite{Delcroix_2025_effective_YTD}.
The following theorem unexpectedly makes the asymptotic problem finite as well.

\begin{thm}
\label{main_thm_rk1}
Let \((X,L)\) be a K-polystable polarized spherical \(G\)-variety of rank one. 
Asymptotic Chow semistability and polystability can each be decided by finitely many explicit polynomial sign and vanishing tests.
\end{thm}

The proof reduces the problem to breakpoints approaching the endpoints of the moment interval. 
Explicit endpoint polynomials control their Chow weights, and higher-order terms are needed only at the (finitely many) positive integer zeros of these polynomials. 

Our rank-one criterion reveals contrasting phenomena. 
For a smooth wonderful \(G_2\)-sevenfold, we show that K-polystability and asymptotic Chow polystability are equivalent for ample polarizations. 
Nevertheless, K-polystability does not imply Chow semistability at a fixed level, and there is no uniform exponent beyond which all K-polystable polarizations become Chow polystable. 
On the other hand, the singular non-toric horospherical fourfold \(Q^4/\mu_6\) provides an example of a K-polystable variety for which the Chow weights of all \(G\)-equivariant product test configurations vanish at every level, but which is not asymptotically Chow semistable.
In this case, instability is detected by a higher-order term in the endpoint expansion, illustrating the necessity of the additional conditions in our criterion.

To conclude this introduction, let us outline the organization of the article. 
Sections~\ref{sec:spherical_var} and~\ref{sec:Chow} recall the spherical data used throughout the paper and explain how Chow stability reduces to the action of the centralizer of \(G\). 
In Section~\ref{sec:spherical_secondary} we construct the spherical secondary polytope, and in Section~\ref{sec:equiv_Chow_polytope} we identify it with the equivariant Chow polytope. 
Theorem~\ref{main_thm_Chow_polytope} is proved there and these two sections constitute the heart of the first part of the article. 
Section~\ref{sec:CSSV} derives the convex criterion for Chow stability of Theorem~\ref{main_thm_Chow_stability} and relates its asymptotics to K-stability.
Section~\ref{sec:rk1} is devoted to spherical rank one and contains the precise statement and proof of Theorem~\ref{main_thm_rk1}. 
Finally, Section~\ref{sec:degen} describes the components of the limit cycles in terms of satellites and uses them to compute exposed faces of the equivariant Chow polytope.

\section{Spherical varieties}
\label{sec:spherical_var}

We recall the facts about spherical varieties needed for the equivariant Chow polytope computation. 
For general background on spherical varieties, we refer to Perrin's survey~\cite{Perrin_2014}.

\subsection{Definition}

Let \(G\) denote a complex connected reductive group. 
Let \(B\subset G\) denote a Borel subgroup, and let \(T\subset B\) denote a maximal torus. 
We denote by \(X^*(T)\) the group of characters of \(T\) (and similarly, \(X^*(B)\) for the group of characters of \(B\), etc.). 
We denote by \(X_*(T)\) the group of one-parameter subgroups, and will in general denote the duality by a bracket \(\langle \cdot,\cdot \rangle\). 

\begin{defn}
A normal algebraic \(G\)-variety \(X\) is a \emph{spherical variety} if it contains an open \(B\)-orbit. 
\end{defn}

A normal algebraic \(G\)-variety with an open \(B\)-orbit also contains an open \(G\)-orbit. 
Choosing a point in this open orbit identifies it with \(G/H\), 
called a \emph{spherical homogeneous space}; the stabilizer \(H\) is called a \emph{spherical subgroup} of \(G\).
 
We will focus on projective spherical varieties, that is on \(G\)-equivariant subvarieties \(X\subset \mathbb{P}(V)\), where \(V\) is a finite-dimensional \(G\)-representation, and such that \(X\) is spherical as a \(G\)-variety. 
Most of the time, we will focus further on the case in which this embedding is the Kodaira embedding for a very ample, \(G\)-linearized line bundle on \(X\). 
We recall that given a Cartier divisor \(\sum n_D D\) on \(X\), up to passing to a multiple, there exists a \(G\)-linearization of the associated line bundle~\cite{KKV_1989}, and different possible linearizations differ by a character of \(G\). 

On a \(G\)-variety \(X\), the action of \(G\) on functions is by \((g\cdot f)(x)=f(g^{-1}\cdot x)\). 
For a \(G\)-linearized line bundle on \(X\), there is an induced action of \(G\) on global sections by \((g\cdot s)(x) = g\cdot s (g^{-1}\cdot x)\), consistent with seeing functions as sections of the trivial bundle with trivial linearization. 

\subsection{Divisors, line bundles and polytopes}

Let \(X\) be a \(G\)-spherical variety.  
Let \(\mathcal{D}(X)\) denote its set of \emph{colors}, which are the \(B\)-stable but not \(G\)-stable prime divisors, and let \(\mathcal{B}(X)\) denote the set of \(G\)-stable prime divisors in \(X\).
Note that the colors of \(X\) are the closures in \(X\) of the colors of the open homogeneous space \(G/H\).

Let \(M\) denote the \emph{weight lattice} of \(X\), that is, the set of  \(B\)-weights of \(B\)-eigenfunctions in \(\mathbb{C}(X)\), regarded as a subgroup of \(X^*(B)=X^*(T)\). 
We also set \(N=\operatorname{Hom}(M,\mathbb{Z})\), 
\(M_\mathbb{R}=M\otimes_\mathbb{Z}\mathbb{R}\)  and \(N_\mathbb{R}=N\otimes_\mathbb{Z}\mathbb{R}\).
For \(m\in M\), let \(f_m\) denote a non-zero \(B\)-eigenfunction of weight \(m\) so that 
\(f_m(b\cdot x)=m(b)^{-1}f_m(x)\). 
Such a function is uniquely defined up to a non-zero constant since \(B\) acts with an open orbit. 
Given a prime \(B\)-stable divisor \(D\), define an element \(\nu_D\in N\) by \(\langle \nu_D,m\rangle:=\operatorname{ord}_D(f_m)\).

\begin{defn}
Let \(X\) be a spherical \(G\)-variety.
The spherical rank of \(X\), denoted by \(\operatorname{rank}(X)\), is the rank of the abelian group \(M\). 
\end{defn}

When no confusion is possible, in this article we shall simply denote the rank of \(X\) by \(r\).
Given a \(B\)-stable Cartier divisor 
\(\delta = \sum_{D\in \mathcal{D}(X)\cup \mathcal{B}(X)} n_D D \), 
we define the associated (possibly empty) polyhedron 
\[ P_{\delta} := \left\{ m\in M_{\mathbb{R}} \ \middle|\ \langle\nu_D,m\rangle+n_D\geq 0 \text{ for all } D\in\mathcal{D}(X)\cup\mathcal{B}(X) \right\}. \]

\begin{thm}[\cite{Brion_1989}]
Let \(L\) be a \(G\)-linearized line bundle on \(X\) and let \(s\) be a rational \(B\)-eigensection with \(B\)-weight \(\chi\). 
Write \(\operatorname{div}(s) = \sum_{D\in\mathcal{D}(X)\cup\mathcal{B}(X)}n_DD\).
Then 
\[ H^0(X,L^{\otimes k}) = \bigoplus_{m\in kP_{\operatorname{div}(s)}\cap M} W_{k,k\chi+m}, \]
where \(W_{k,\mu}\) denotes the irreducible submodule of highest-weight \(\mu\) occurring in \(H^0(X,L^{\otimes k})\).
\end{thm}

We define the moment polytope of the polarized \(G\)-variety \((X,L)\) by
\[
\Delta(X,L):=\chi+P_{\operatorname{div}(s)}.
\]
This polytope depends only on the \(G\)-linearized line bundle \(L\), and not on a choice of section. 
Thus the set of highest-weights occurring in \(H^0(X,L^{\otimes k})\) is
\[
\Lambda_k
:=
k\Delta(X,L)\cap(k\chi+M),
\]
and the preceding decomposition may equivalently be written
\[
H^0(X,L^{\otimes k})
=
\bigoplus_{\lambda\in\Lambda_k}W_{k,\lambda}.
\]
When \(L\) is ample, the dimension of the moment polytope coincides with the spherical rank of \(X\).

\subsection{The valuation cone and spherical roots}

Let \(G/H\subset X\) be the open \(G\)-orbit.
A \(G\)-invariant discrete valuation of \(\mathbb{C}(G/H)\),
trivial on \(\mathbb{C}^*\), determines an element of \(N_\mathbb{Q}\)
by restriction to \(B\)-eigenfunctions.
The cone generated by these elements is the
\emph{valuation cone}
\(
\mathcal V\subset N_\mathbb{R}.
\)
This is a co-simplicial cone: there exists a finite set of linearly independent primitive elements
\(
\Sigma\subset M,
\)
called the \emph{spherical roots} of \(G/H\), such that
\[
\mathcal V
=
\left\{
v\in N_\mathbb{R}
\;\middle|\;
\langle\alpha,v\rangle\leq0
\text{ for every }\alpha\in\Sigma
\right\}.
\]
Throughout the paper we use the negative-polar convention
\[
\mathcal V^\vee
=
\left\{
\gamma\in M_\mathbb{R}
\;\middle|\;
\langle\gamma,v\rangle\leq0
\text{ for every }v\in\mathcal V
\right\},
\]
so that
\(
\mathcal V^\vee=\operatorname{cone}(\Sigma).
\)
We will also denote by \(\operatorname{Lin}(\mathcal{V})\) the linear part of the valuation cone, defined by 
\[ \operatorname{Lin}(\mathcal{V}) = \mathcal{V}\cap -\mathcal{V}. \]

Let \(L\) be a \(G\)-linearized ample line bundle, the multiplication in the section ring 
\(
H^0(X,L^{\otimes k})
=
\bigoplus_{\lambda\in\Lambda_k}W_{k,\lambda}
\)
 satisfies
\begin{equation}
\label{eq:section-ring-multiplication}
W_{k,\lambda}W_{\ell,\mu}
\subseteq
\bigoplus_{\substack{
\gamma\in\mathcal V^\vee\\
\lambda+\mu-\gamma\in\Lambda_{k+\ell}}}
W_{k+\ell,\lambda+\mu-\gamma}.
\end{equation}
We call \(\gamma\) a \emph{tail} of the corresponding
product component.
This is the section-ring form of the characterization of the valuation cone by inequalities arising from multiplication, see~\cite[Lemma~6.1]{Knop_1991}.

\begin{rem}
We call \(X\) \emph{horospherical} if \(\mathcal V=N_\mathbb{R}\), equivalently 
\(
\Sigma=\varnothing\) or
\( \mathcal V^\vee=\{0\}\).
In this case, multiplication of sections contains only the component corresponding to the sum of weights, called the \emph{Cartan component}: \(W_{k,\lambda}W_{\ell,\mu}
=
W_{k+\ell,\lambda+\mu}\).
In the general spherical case, the multiplication of sections always yields a non-zero contribution in the Cartan component. 
\end{rem}

\subsection{The Duistermaat--Heckman polynomial}

We recall the degree formula for polarized spherical varieties due to Brion~\cite[ Th\'eor\`eme~4.1]{Brion_1989}.
Let \(\Phi\) be the root system of \((G,T)\),
let \(\Phi^+\) be the set of positive roots
determined by \(B\), and let \(\rho\) denote
the half-sum of the positive roots.
Set
\[
E=
\left\{
\alpha\in\Phi^+ 
\;\middle|\;
\langle x,\alpha^\vee\rangle=0
\text{ for every }x\in\Delta(X,L)
\right\}
\]
and define the Duistermaat--Heckman polynomial as
\[
P_{\mathrm{DH}}(x)
=
\prod_{\alpha\in\Phi^+\setminus E}
\frac{\langle x,\alpha^\vee\rangle}
     {\langle\rho,\alpha^\vee\rangle}
\]
for \(x\in X^*(T)\). 
The Weyl dimension formula gives, uniformly in 
\(\lambda/k\in\Delta(X,L)\),
\[
\dim W_{k,\lambda}
=
k^{n-r}
P_{\mathrm{DH}}\left(\frac{\lambda}{k}\right)
+
O(k^{n-r-1}).
\]
We denote by \(d\mu_M\) the Lebesgue measure on the affine span of \(\Delta(X,L)\), normalized by the lattice \(M\) in the direction of \(\Delta(X,L)\).
Brion's formula expresses the degree of the ample line bundle \(L\) in terms of this data: 
\[
L^n
=
n!\int_\Delta P_{\mathrm{DH}}(x)\,d\mu_M.
\]

\section{Chow stability and equivariant Chow polytope}
\label{sec:Chow}

In this section, we recall the construction of the Chow form and the definition of Chow stability~\cite{Mumford_1977}. 
For a \(G\)-stable subvariety of the projectivization of a multiplicity-free \(G\)-module, Kempf's instability theorem allows us to test Chow semistability and polystability using the connected centralizer in \(\operatorname{SL}(V)\). 
Following the weight-polytope construction of~\cite{Kapranov_Sturmfels_Zelevinsky_1992}, we also consider the Chow weights for the centralizer in \(\operatorname{GL}(V)\), which acts by independent scalars on the irreducible summands.

\subsection{The Chow point and Chow stability}

Let \(V\) be a complex vector space of dimension \(N+1\), and let \(X\subset\mathbb P(V)\) be an irreducible projective subvariety of dimension \(n<N\) and degree \(d\).
Let \(\operatorname{Gr}(N-n,V)\) denote the Grassmannian of \((N-n)\)-dimensional vector subspaces of \(V\). 
The \emph{associated hypersurface} of \(X\) is
\[
Z_X = \{L\in\operatorname{Gr}(N-n,V) \mid X\cap\mathbb P(L)\ne \emptyset\}.
\]

We consider the dual Pl\"ucker embedding of the Grassmannian
\[
\operatorname{Gr}(N-n,V)
\simeq \operatorname{Gr}(n+1,V^*)
\hookrightarrow \mathbb P(\bigwedge^{n+1}V^*).
\]
The associated hypersurface of \(X\) is defined by a section
\[
R_X\in
\mathcal B_d:=
H^0(\operatorname{Gr}(n+1,V^*),\mathcal O(d)),
\]
unique up to a nonzero scalar, called the
\emph{Chow form} of \(X\).
The point \([R_X]\in\mathbb P(\mathcal B_d)\) is called the \emph{Chow point} of \(X\).

The Chow form can equivalently be viewed as a polynomial in \(n+1\) linear forms on \(V\), by setting
\[
R_X(f_1,\ldots,f_{n+1})
=
R_X(f_1\wedge\cdots\wedge f_{n+1})
\]
where \(f_1\wedge\cdots\wedge f_{n+1}\) denotes the \((n+1)\)-dimensional subspace of \(V^*\) spanned by \(f_1,\ldots,f_{n+1}\). 
This polynomial is called the \emph{\(X\)-resultant}, or equivalently the Chow form in resultant notation.
It is homogeneous of degree \(d\) in each argument, and satisfies
\[
R_X(f_1,\ldots,f_{n+1})=0
\quad\Longleftrightarrow\quad
X\cap\{f_1=\cdots=f_{n+1}=0\}
\neq\varnothing.
\]

The natural action of \(\operatorname{GL}(V)\) on \(V^*\) induces a canonical linearization of the Pl\"ucker line bundle \(\mathcal O(1)\) on \(\operatorname{Gr}(n+1,V^*)\), and hence a linear representation on \(\mathcal B_d\).
Restricting this representation to \(\operatorname{SL}(V)\), we obtain the linearization used to define Chow stability.

\begin{defn}
The projective subvariety \(X\subset \mathbb P(V)\) is Chow unstable (resp. semistable,  polystable) if \(R_X\in \mathcal{B}_d\) is GIT unstable (resp. semistable, polystable) with respect to the action of \(\operatorname{SL}(V)\). 
\end{defn}

\begin{rem}
Recall the general definition of GIT stability: given a reductive group \(G\) and a finite dimensional representation \(V\) of \(G\), a point \([v]\in \mathbb P(V)\) is: 
\begin{itemize}
\item unstable if \(0\) is in the closure of the orbit \(G\cdot v\), 
\item semistable if \(0\) is not in the closure of \(G\cdot v\), 
\item polystable if the orbit \(G\cdot v\) is closed.
\end{itemize}
\end{rem}

More generally, if \(Z=\sum_i m_i[X_i]\) is an effective cycle of pure dimension \(n\), we define its Chow form by
\[
R_Z=\prod_i R_{X_i}^{m_i}.
\]
The Chow stability notions extend to effective cycles through their Chow points.

\subsection{Kempf's instability and the equivariant Chow polytope}

Assume that \(V\) is a \emph{multiplicity-free} representation of a connected reductive group \(G\),  that is \(V = \bigoplus_{i\in I} V_i\), where the \(V_i\) are pairwise non-isomorphic irreducible \(G\)-representations. 

By Schur's lemma, the centralizer \(
\widetilde T_0:=Z_{\operatorname{GL}(V)}(G)
\) is a torus isomorphic to \((\mathbb{C}^*)^I\).
Its intersection with \(\operatorname{SL}(V)\) has identity component
\( T_0:=Z_{\operatorname{SL}(V)}(G)^\circ\),
a torus of dimension \(|I|-1\).

\begin{prop}\label{prop:Chow-centralizer}
Let \(X\subset\mathbb P(V)\) be a \(G\)-stable projective subvariety, where \(V\) is a multiplicity-free \(G\)-module.
Then \(X\) is Chow semistable, respectively Chow polystable, if and only if its Chow point is semistable, respectively polystable, for the action of \(T_0\).
\end{prop}

\begin{proof}
Let \(G'\) be the identity component of the inverse image of the projective image of \(G\) under \(\operatorname{SL}(V)\to\operatorname{PGL}(V)\).
Then \(G'\) is reductive and fixes the Chow point \([R_X]\).
The images of \(G\) and \(G'\) in \(\operatorname{PGL}(V)\) agree, so their lifts differ only by scalars and have the same centralizer in \(\operatorname{SL}(V)\).

We apply the centralizer criterion from Kempf's instability theorem~\cite[Corollary~4.5]{Kempf_1978}, in the formulation for semistability and polystability given in~\cite[Theorem~2.2 and Corollary~2.3]{Lee_Li_Sturm_Wang_2019}.
Both assertions thus reduce to the action of \(Z_{\operatorname{SL}(V)}(G')\).
Passing to its identity component \(T_0\) preserves semistability and polystability, since the full centralizer orbit of \(R_X\) is a finite union of translates of its \(T_0\)-orbit.
\end{proof}

Write the weight decomposition of the Chow form \(R_X\) in the \(T_0\)-representation \(\mathcal{B}_d\) as
\[
R_X=\sum_{\chi\in X^*(T_0)}R_{X,\chi}.
\]

\begin{defn}
The \emph{equivariant Chow polytope} \(\operatorname{Chow}_{\mathrm{SL}}^G(X)\) is the convex hull of the \(T_0\)-weights occurring in the weight decomposition of the Chow form \(R_X\in\mathcal B_d\): 
\[
\operatorname{Chow}_{\mathrm{SL}}^G(X)
=
\operatorname{conv}\left\{
\chi\in X^*(T_0)\mid R_{X,\chi}\ne0
\right\}
\subset X^*(T_0)\otimes_{\mathbb Z}\mathbb R.
\]
\end{defn}

The equivariant Chow polytope is also the moment polytope of the normalization of the projective torus orbit closure \( \overline{T_0\cdot[R_X]} \subset\mathbb P(\mathcal B_d)\), equipped with the pullback of \(\mathcal O(1)\) and its induced \(T_0\)-linearization.
In general, this construction is called the weight polytope of \(R_X\) with respect to the action of \(T_0\). 

By Proposition~\ref{prop:Chow-centralizer}, it suffices to consider the action of \(T_0\). 
The Hilbert--Mumford criterion~\cite[\S~1.1]{Mumford_1977} gives semistability when the weight polytope contains the origin. 
In the semistable case, a one-parameter subgroup with minimum weight zero has an affine limit given by retaining the weight components on the corresponding face~\cite[Proposition~1.4]{Kapranov_Sturmfels_Zelevinsky_1992}.
Such a limit leaves the orbit precisely when that face is proper. 
Thus polystability requires the origin to lie in the relative interior (denoted by \(\operatorname{relint}\), the interior
in the affine span), and we obtain
\[
\begin{aligned}
X\text{ is Chow semistable}
&\iff 0\in\operatorname{Chow}_{\mathrm{SL}}^G(X),\\
X\text{ is Chow polystable}
&\iff 0\in\operatorname{relint}\operatorname{Chow}_{\mathrm{SL}}^G(X).
\end{aligned}
\]

In this paper, it will be convenient to work as well with the weight polytope of \(R_X\) with respect to the action of the larger torus \(\widetilde{T}_0\) in \(\mathrm{GL}(V)\), which we denote by 
\[ \operatorname{Chow}_{\mathrm{GL}}^G(X) \subset X^*(\widetilde{T}_0)\otimes\mathbb{R}. \]

\begin{exa}\label{exa:flag-variety-Chow}
Let \(X=G/P\) be a positive-dimensional generalized flag variety, and let \(L\) be a very ample \(G\)-linearized line bundle. 
By the Borel--Weil theorem, the representation \(V=H^0(X,L)^*\) is irreducible. 
Hence \(\widetilde T_0\simeq\mathbb G_m\) acts by scalars and \(T_0\) is trivial.
Writing \(\varepsilon\) for the character of scalar multiplication on \(V\), and \(n=\dim X\), we obtain
\( \operatorname{Chow}_{\mathrm{GL}}^G(X)
=\{(n+1)L^n\,\varepsilon\}\) and
\( \operatorname{Chow}_{\mathrm{SL}}^G(X)=\{0\}\).
The complete embedding defined by \(L\) is therefore Chow polystable (but not Chow stable, since its projective stabilizer contains the positive-dimensional image of \(G\)).
\end{exa}

\section{The spherical secondary polytope}
\label{sec:spherical_secondary}

\subsection{Decomposition polytopes and spherical envelopes}

Let \(X\) be a projective spherical \(G\)-variety of dimension \(n\), and let \(L\) be an ample \(G\)-linearized line bundle on \(X\).
We write \(\Delta=\Delta(X,L)\) for the moment polytope and \(P_{\mathrm{DH}}\,d\mu_M\) for the Duistermaat--Heckman measure.
The tail cone is \(\mathcal V^\vee=\operatorname{cone}(\Sigma)\).

The following characterization of base-point-free \(G\)-stable linear systems for polarized spherical varieties follows directly from the proof of~\cite[Theorem~4.14]{Cupit-Foutou_Pezzini_Van_Steirteghem_2020}. 

\begin{prop}\label{prop:decomposition-polytopes-nonempty}
Let \( E=\bigoplus_{\lambda\in J}W_{1,\lambda} \) be a non-zero \(G\)-stable subspace of \(H^0(X,L)\). 
It is base-point-free if and only if 
\( \Delta\subset\operatorname{conv}(J)-\mathcal V^\vee\).
In this case, the morphism defined by \(E\) is finite onto its image. 
\end{prop}

We now fix a base-point-free \(G\)-invariant linear system \(E\subset H^0(X,L)\) whose associated morphism
\[
\pi:X\longrightarrow Y\subset\mathbb P(E^*)
\]
is finite onto its image. 
Thus \(L=\pi^*\mathcal O_Y(1)\). 
The pair \((X,\pi)\) is a spherical variety over \(\mathbb P(E^*)\) in the sense of~\cite[\S~1.2]{Alexeev_Brion_2006}.
The linear system \(E\) is part of the data and need not be complete in this section.
We write 
\begin{equation}
\label{eqn:defn_R_A_E}
R=\bigoplus_{k\geq0}H^0(X,L^{\otimes k}),
\qquad
A_E=\operatorname{im}\!\left(\operatorname{Sym}E\longrightarrow R\right)
\end{equation}
for the section ring of \((X,L)\) and the homogeneous coordinate
ring of \(Y\subset\mathbb P(E^*)\), respectively.

As in the statement above, we write \( E=\bigoplus_{\lambda\in J}W_{1,\lambda} \), where \(J\subset\Lambda_1\) is the set of highest-weights occurring in \(E\).
Since \(E\) is multiplicity-free, the centralizer \( \widetilde T_0=Z_{\mathrm{GL}(E^*)}(G) \) is the torus of independent scalar multiplications on the irreducible summands \(W_{1,\lambda}^*\).
We use \(\mathbb R^J\) as the full block-coordinate space.

\begin{defn}
\label{def:decomposition-polytope}
For \(x\in\Delta\), define the \emph{decomposition polytope at \(x\)} as
\[
K_x=
\left\{
r\in\mathbb R_{\geq0}^{J}
\;\middle|\;
\sum_{\lambda\in J}r_\lambda=1,\quad
\sum_{\lambda\in J}r_\lambda\lambda-x\in \mathcal V^\vee
\right\}.
\]
\end{defn}

This encodes the various ways to write \(x\) as a convex combination of elements of the weights in \(J\) up to an element of the tail cone \(\mathcal V^\vee\). 
It follows directly from the preceding proposition that \(K_x\) is nonempty for every \(x\in\Delta\).

\begin{defn}
The \emph{support function} \(s_Q:\mathbb{R}^J\to \mathbb{R}\) of a nonempty compact convex set \(Q\subset\mathbb R^J\) is defined by 
\[
s_Q(h)=\max_{u\in Q}\langle u,h\rangle,
\]
for \(h\in\mathbb R^J\),  where
\(\langle u,h\rangle=\sum_{\lambda\in J}u_\lambda h_\lambda\).
\end{defn}

The function \(x\mapsto s_{K_x}(h)\) for a fixed height vector \(h\), where \(s_{K_x}\) is the support function of the decomposition polytope, will play an important role in the article.

\begin{defn}
For \(h=(h_\lambda)_{\lambda\in J}\in\mathbb R^J\), the \emph{spherical envelope} is the function \(g_h\) on \(\Delta\) defined by 
\[
g_h(x):=s_{K_x}(h)
=
\max_{r\in K_x}\sum_{\lambda\in J}r_\lambda h_\lambda,
\qquad \text{for } x\in\Delta.
\]
\end{defn}

\begin{lem}\label{lem:spherical-envelope-duality}
For every \(x\in\Delta\),
\[
g_h(x)
=
\min_{\substack{
b\in\mathbb R,\ w\in\mathcal V\\
b+\langle\lambda-\chi,w\rangle\geq h_\lambda
\quad(\lambda\in J)}}
\bigl(b+\langle x-\chi,w\rangle\bigr).
\]
In particular, \(g_h\) is concave and piecewise-affine,
with slopes in \(\mathcal V\).
If \(h\in\mathbb Q^J\), its affine pieces are rational.
\end{lem}

\begin{proof}
Set
\[
P_h=\operatorname{conv}\{(\lambda-\chi,h_\lambda) \mid \lambda\in J\}
    +\bigl(-\mathcal V^\vee\times\mathbb R_{\leq0}\bigr).
\]
By the definition of \(K_x\), for every \(x\in\Delta\),
\[
g_h(x)=\max\{y \mid (x-\chi,y)\in P_h\}.
\]
The set \(P_h\) is a polyhedron~\cite[Theorem~19.1]{Rockafellar_1970}. 
Since it is closed under translation in the negative vertical direction, its defining inequalities are either independent of \(y\) or of the form \(y\leq b+\langle z,w\rangle\).

An inequality of the latter form is valid on \(P_h\) if and only if it holds at each \((\lambda-\chi,h_\lambda)\) and \(\langle\gamma,w\rangle\leq0\) for every \(\gamma\in\mathcal V^\vee\), equivalently \(w\in\mathcal V\).
Thus \(g_h\) is the pointwise minimum of the affine functions in the statement, and finitely many suffice. 
In particular, \(g_h\) is concave and piecewise affine. 
If \(h\) is rational, then \(P_h\) is a rational polyhedron, so the affine functions can be chosen with rational coefficients.
\end{proof}

This description explains the term envelope: in the horospherical case \(\mathcal{V}=N\otimes \mathbb{R}\), the spherical envelope is the smallest concave function whose values at \(\lambda\in \Delta\) are larger than or equal to \(h_{\lambda}\). 

\subsection{The spherical secondary polytope}

To define the spherical secondary polytope, we will use the notion of Minkowski integral, which generalizes Minkowski sums of convex sets. 

\begin{defn}
Let \((S,\nu)\) be a finite measure space and \(V\) a finite-dimensional real vector space.
Given a measurable, integrably bounded family \((Q_x)_{x\in S}\) of nonempty compact convex subsets of \(V\), its \emph{Minkowski integral} 
\(
\int_S Q_x\,d\nu(x)
\)
is the unique compact convex subset of \(V\) whose support function is given, for \(h\in V^*\) by
\[
s_{\int_S Q_x\,d\nu(x)}(h)
=\int_S s_{Q_x}(h)\,d\nu(x).
\]
\end{defn}

Equivalently, the Minkowski integral is the set of integrals of measurable selections of the family \((Q_x)_{x\in S}\), also known as the Aumann integral~\cite{Aumann_1965}.

\begin{defn}
The \emph{spherical secondary polytope} associated with \((X,L,E)\) is
\[
\Sigma_{\mathrm{sec}}^G(X,L;E)
=
(n+1)!\int_\Delta K_x\,P_{\mathrm{DH}}(x)\,d\mu_M(x)
\subset\mathbb R^J.
\]
For the complete linear system \(E=H^0(X,L)\), we write \(\Sigma_{\mathrm{sec}}^G(X,L)\).
\end{defn}

\begin{prop}\label{prop:spherical-secondary-polytope}
The convex set \(\Sigma_{\mathrm{sec}}^G(X,L;E)\) is a polytope.
\end{prop}

\begin{proof}
Choose a finite set of generators \(v_1,\ldots,v_s\) of \(\mathcal V\) as a cone, including both signs of a basis of its linear part. 
Using the negative-polar convention for \(\mathcal V^\vee\), we may write
\[
K_x=
\left\{r\in\mathbb R^J:
Ar\leq b(x),\quad \sum_{\lambda\in J}r_\lambda=1
\right\},
\]
where the inequalities \(Ar\leq b(x)\) are
\[
r_\lambda\geq0\quad(\lambda\in J),
\qquad
\sum_{\lambda\in J}\langle\lambda-\chi,v_i\rangle r_\lambda
\leq \langle x-\chi,v_i\rangle
\quad(1\leq i\leq s).
\]
Thus \(A\) is independent of \(x\).

The normal cone of any nonempty face of \(K_x\) is generated by the rows of \(A\) corresponding to the inequalities active on that face, together with the line spanned by \((1,\ldots,1)\).
Thus only finitely many normal cones can occur as \(x\) varies in \(\Delta\).
Consequently, only finitely many normal fans of the polytopes \(K_x\) can occur.

Let \(\mathcal F\) be a common refinement of these finitely many normal fans.
For every \(x\in\Delta\), the support function \(s_{K_x}\) is linear on each cone of \(\mathcal F\).
By the definition of the Minkowski integral,
\[
\begin{aligned}
s_{\Sigma_{\mathrm{sec}}^G(X,L;E)}(h)
&=(n+1)!\int_\Delta
  s_{K_x}(h)\,P_{\mathrm{DH}}(x)\,d\mu_M(x)\\
&=(n+1)!\int_\Delta
  \max_{r\in K_x}\langle r,h\rangle\,
  P_{\mathrm{DH}}(x)\,d\mu_M(x).
\end{aligned}
\]
This function is therefore linear on each cone of \(\mathcal F\).

Hence the support function of \(\Sigma_{\mathrm{sec}}^G(X,L;E)\) is piecewise-linear on a finite complete polyhedral fan.
Therefore, the corresponding compact convex body is a polytope; see~\cite[Section~1.7]{Schneider_2014}.
\end{proof}

\begin{rem}[The spherical secondary polytope as a weighted fiber polytope]
\label{rem:BS}
The polytope
\[
P=\{(x,r)\in\Delta\times\mathbb R^J\mid r\in K_x\}
\]
projects onto \(\Delta\) with fibers \(K_x\), via the identification \((x,r)\mapsto r\) on each fiber.
The ordinary fiber polytope of this projection~\cite{Billera_Sturmfels_1992} has the same normal fan as \(\Sigma_{\mathrm{sec}}^G(X,L;E)\).
More precisely, the spherical secondary polytope may be viewed, up to normalization, as a weighted fiber polytope, with respect to the Duistermaat--Heckman measure, in the sense of
\cite{Forcier_Gaubert_Leclere_2024}. 
This also provides a way to compute the spherical secondary polytope, following the  arguments of~\cite{Billera_Sturmfels_1992}. 

Let \(\pi:P\to\Delta\) denote the projection.
For \(x\in\Delta\), define
\[
\sigma_x=
\bigcap_{\substack{F\text{ face of }P\\ x\in\pi(F)}}
\pi(F).
\]
The chamber complex \(\mathcal C\) of the projection \(\pi\), in the sense of the previous references, is the polyhedral subdivision of \(\Delta\) formed by the polytopes \(\sigma_x\), for \(x\in\Delta\).
Equivalently, it is the coarsest polyhedral subdivision of \(\Delta\) for which the image of every face of \(P\) is a union of cells.

The chamber complex \(\mathcal C\) can equivalently be described as the common refinement of the subdivisions of \(\Delta\) into maximal domains of affinity of the spherical envelopes \(g_h\), for all \(h\in\mathbb R^J\).
Indeed, every face \(F\) of \(P\) can be written as \(F=\{(x,r)\in P\mid \langle h,r\rangle=\ell(x)\}\), where \(h\in\mathbb R^J\) and \(\ell\) is an affine function satisfying
\(\langle h,r\rangle\leq\ell(x)\) for all \((x,r)\in P\).
Since \(g_h(x)=\max_{r\in K_x}\langle h,r\rangle\), we have
\(
\pi(F)=\{x\in\Delta\mid g_h(x)=\ell(x)\}.
\)
Thus the image of every face of \(P\) is a cell of the subdivision induced by some spherical envelope.
Conversely, every maximal linearity domain of a spherical envelope \(g_h\) is obtained in this way, by taking for \(\ell\) the affine restriction of \(g_h\) on this domain. 
Note that the chamber complex \(\mathcal C\) is independent of the Duistermaat--Heckman measure and need not itself be induced by a spherical envelope.

For each maximal cell \(\sigma\in\mathcal C\), let
\[
m_\sigma=\int_\sigma P_{\mathrm{DH}}(x)\,d\mu_M(x),
\qquad
b_\sigma=\frac{1}{m_\sigma}
\int_\sigma xP_{\mathrm{DH}}(x)\,d\mu_M(x)
\]
denote its Duistermaat--Heckman volume and barycenter, respectively.
Then
\[
\Sigma_{\mathrm{sec}}^G(X,L;E)
=(n+1)!\sum_{\sigma\in\mathcal C^{\max}}
m_\sigma K_{b_\sigma},
\]
where the sum is a Minkowski sum.
Indeed, since \(g_h\) is affine on each cell \(\sigma\), we have
\[
\int_\sigma g_h(x)P_{\mathrm{DH}}(x)\,d\mu_M(x)
=m_\sigma g_h(b_\sigma).
\]
The formula follows by identifying the support functions of both sides.
\end{rem}

\subsection{Examples of spherical secondary polytopes}

\begin{exa}[The toric conic]\label{exa:toric-conic}
Let \(X=\mathbb P^1\), \(L=\mathcal O(2)\), with the standard \(\mathbb G_m\)-action. For the complete linear system, we identify \(\Delta=[0,2]\) and \(J=\{0,1,2\}\).
The valuation and tail cone are trivial, and \(P_{\mathrm{DH}}\equiv 1\). 
For \(x\in[0,2]\), the decomposition polytope is
\[
K_x=\left\{
(1-x+t,x-2t,t):
\max(0,x-1)\leq t\leq x/2
\right\},
\]
a segment whose two endpoints depend piecewise affinely on \(x\).
Integrating these gives
\[
\Sigma_{\mathrm{sec}}^G(X,L)
=2\int_0^2 K_x\,dx
=\operatorname{conv}\{(2,0,2),(1,2,1)\}.
\]
This is the classical secondary polytope of the three-point configuration \(\{0,1,2\}\).

The secondary polytope highly depends on the choice of linear system \(E\): now consider the smaller linear system \(E=\langle u^2,v^2\rangle\).
It defines the degree-two morphism \(\mathbb P^1\longrightarrow\mathbb P^1\), \([u:v]\longmapsto[u^2:v^2]\).
Here \(J=\{0,2\}\), and every decomposition polytope is a singleton \(K_x=\left\{\left(1-\frac{x}{2},\frac{x}{2}\right)\right\}\). Consequently, \(\Sigma_{\mathrm{sec}}^G(X,L;E)=\{(2,2)\}\).
\end{exa}

\begin{exa}[A horospherical projective space]
\label{exa:horospherical-P3}
Let \(G=\mathrm{SL}_2\times\mathrm{SL}_2\) act on \(X=\mathbb P(\mathbb C^2\oplus\mathbb C^2)\) polarized by \(L=\mathcal O_X(2)\).
This is a rank-one horospherical threefold.  
Denote by \(\omega_1,\omega_2\) the fundamental weights
of the two factors. The spherical lattice is \(M=\mathbb Z(\omega_1-\omega_2)\) and the moment polytope is \(\Delta=[2\omega_2,2\omega_1] \subset X^*(T)_\mathbb R\) (see e.g.~\cite{Delcroix_Montagard}). 
For the complete linear system, \(J=\{2\omega_2,\omega_1+\omega_2,2\omega_1\}\). 
Writing a point \(x\in \Delta\) as \(s\omega_1+(2-s)\omega_2\), the decomposition polytope is
\[
K_x=\left\{
(1-s+t,s-2t,t):
\max(0,s-1)\leq t\leq s/2
\right\}.
\]
Thus \(K_x\) coincides with the decomposition polytope of Example~\ref{exa:toric-conic}.

However, here the Duistermaat--Heckman measure is \(P_{\mathrm{DH}}\,d\mu_M=s(2-s)\,ds\).
Integrating the endpoints of \(K_x\) gives
\[
\Sigma_{\mathrm{sec}}^G(X,L)
=24\int_0^2 K_{s\omega_1+(2-s)\omega_2}\,s(2-s)\,ds
=\operatorname{conv}\{(16,0,16),(6,20,6)\}.
\]
\end{exa}

\begin{exa}[A non-horospherical spherical surface]
\label{exa:diagonal-SL2}
Let \(G=\mathrm{SL}_2\) act diagonally on \(X=\mathbb P^1\times\mathbb P^1\), polarized by
\(L=\mathcal O(1,1)\).
Writing \(\omega\) for the fundamental weight, we have \(M=\mathbb Z(2\omega)\), \(\Sigma=\{2\omega\}\) and \(\Delta=[0,2\omega]\) (see~\cite{Delcroix_Montagard}).
For the complete linear system, \(J=\{0,2\omega\}\).
Writing \(x=2s\omega\in\Delta\), we obtain
\[
K_x=\left\{(1-t,t):s\leq t\leq1\right\}.
\]
The Duistermaat--Heckman measure is \(P_{\mathrm{DH}}\,d\mu_M=2s\,ds\).
Consequently,
\[
\Sigma_{\mathrm{sec}}^G(X,L)
=6\int_0^1 K_{2s\omega}\,2s\,ds
=\operatorname{conv}\{(2,4),(0,6)\}.
\]
\end{exa}

\begin{exa}[Rank-two spherical quadrics]
\label{exa:two-block-quadric}
Let \(U,W\) be vector spaces of dimensions \(p,q\geq3\),
respectively, with nondegenerate quadratic forms \(q_U,q_W\).
Consider the quadric \(X=\{q_U(u)+q_W(w)=0\} \subset\mathbb P(U\oplus W)\), equipped with the polarization \(L=\mathcal O_X(1)\), under the action of \(G=\mathrm{SO}(U)\times\mathrm{SO}(W)\).
This is a spherical variety of rank two (see~\cite{Del_AIF}).

Let \(\lambda_U,\lambda_W\) denote the highest weights of the standard representations of the two factors.
The spherical lattice is \(M=\{i\lambda_U+j\lambda_W: i,j\in\mathbb Z,\ i+j\in2\mathbb Z\}\), the spherical roots are \(2\lambda_U,2\lambda_W\) and the moment polytope is \(\Delta=\operatorname{conv}(0,\lambda_U,\lambda_W)\).
For the complete linear system, \(J=\{\lambda_U,\lambda_W\}\).
Writing \(x=s\lambda_U+t\lambda_W\in\Delta\), we obtain
\[
K_x=\{(r,1-r):s\leq r\leq1-t\}.
\]
Note that the decomposition polytopes can be non-empty even at points of \(\Delta\) outside \(\operatorname{conv}(J)\), thanks to the nontrivial tail cone. 
This is a major difference with respect to the toric (and more generally horospherical) case.

The Duistermaat--Heckman measure is
\[
P_{\mathrm{DH}}\,d\mu_M
=\frac{2s^{p-2}t^{q-2}}{(p-2)!(q-2)!}\,ds\,dt.
\]
Consequently,
\[
\begin{aligned}
\Sigma_{\mathrm{sec}}^G(X,L)
&=(p+q-1)!
\int_{\substack{s,t\geq0\\s+t\leq1}}
K_{s\lambda_U+t\lambda_W}\,
\frac{2s^{p-2}t^{q-2}}{(p-2)!(q-2)!}\,ds\,dt\\
&=\operatorname{conv}
\{(2p-2,2q),(2p,2q-2)\}.
\end{aligned}
\]
\end{exa}

\begin{exa}[The blow-up of \(\mathbb P^2\times\mathbb P^2\) along the diagonal]
\label{exa:blowup-diagonal}
Let \(X\) be the blow-up of \(\mathbb P^2\times\mathbb P^2\) along the diagonal, under the diagonal action of \(G=\mathrm{SL}_3\).
This is a rank-one non-horospherical spherical fourfold (see~\cite{Delcroix_Montagard}).
Writing \(H_1,H_2\) for the pullbacks of the hyperplane classes and \(E\) for the exceptional divisor, we consider the anticanonical polarization \(L=-K_X=3H_1+3H_2-E\).

Let \(\omega_1,\omega_2\) be the fundamental weights of \(\mathrm{SL}_3\), and \(\alpha_2=2\omega_2-\omega_1\).
The spherical lattice is \(M=\mathbb Z\alpha_2\), the spherical root is \(\alpha_2\), and the moment polytope is \(\Delta=[3\omega_1,\omega_1+4\omega_2]\). 
For the complete linear system, \(J=\{\omega_1+4\omega_2,\,
2\omega_1+2\omega_2,\,3\omega_1\}\).
Writing \(x_s=6\omega_2-s\alpha_2\), with \(1\leq s\leq3\), we obtain 
\[
K_{x_s}=
\left\{
(r_1,r_2,r_3)\in\mathbb R_{\geq0}^3:
r_1+r_2+r_3=1,\quad r_1+2r_2+3r_3\leq s
\right\}.
\]
These decomposition polytopes are triangles for \(1<s<2\) and quadrilaterals for \(2<s<3\).

The Duistermaat--Heckman measure is \(P_{\mathrm{DH}}\,d\mu_M=s(3-s)(6-s)\,ds\).
Consequently,
\[
\Sigma_{\mathrm{sec}}^G(X,L)
=120\int_1^3 K_{x_s}\,s(3-s)(6-s)\,ds.
\]
Computing the support function gives the quadrilateral
\[
\Sigma_{\mathrm{sec}}^G(X,L)
=\operatorname{conv}\left\{
\begin{aligned}
&(1680,0,0),\ (606,1074,0),\\
&(606,900,174),\ (1056,0,624)
\end{aligned}
\right\}.
\]
It lies in the affine plane \(p_1+p_2+p_3=1680\).

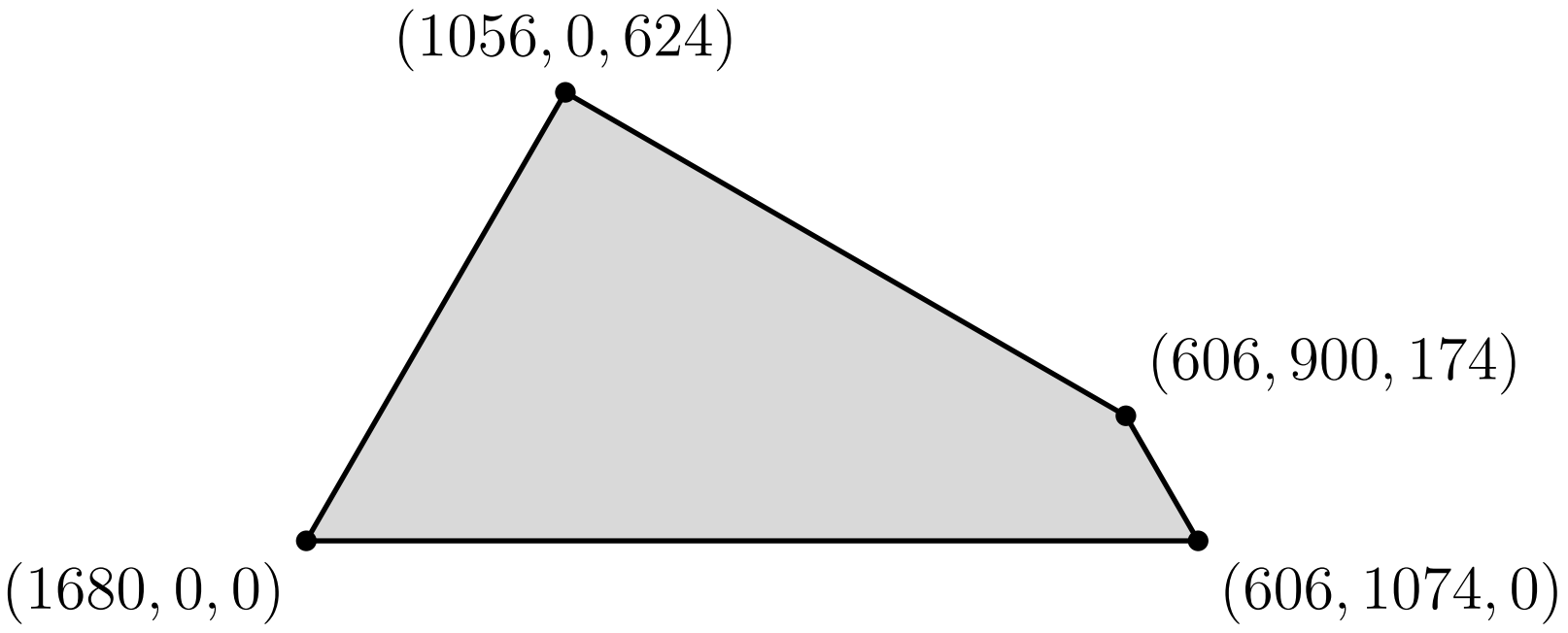
\begin{figure}[ht]
\centering
\begin{tikzpicture}[x=0.004cm,y=0.004cm]
\coordinate (A) at (0,0);
\coordinate (B) at ({1074*sqrt(2)},0);
\coordinate (C) at ({987*sqrt(2)},{174*sqrt(3/2)});
\coordinate (D) at ({312*sqrt(2)},{624*sqrt(3/2)});

\filldraw[
    fill=black!15,
    draw=black,
    line width=1pt
]
(A) -- (B) -- (C) -- (D) -- cycle;

\foreach \v in {A,B,C,D}
    \fill[black] (\v) circle[radius=2pt];
\node[below left]  at (A) {$(1680,0,0)$};
\node[below right] at (B) {$(606,1074,0)$};
\node[above right] at (C) {$(606,900,174)$};
\node[above]       at (D) {$(1056,0,624)$};

\end{tikzpicture}
\caption{The spherical secondary polytope of $\operatorname{Bl}_{\operatorname{Diag}(\mathbb P^2)}(\mathbb P^2\times\mathbb P^2)$ with its anticanonical polarization.}
\label{fig:secondary-blowup-diagonal}
\end{figure}
\end{exa}

\begin{exa}[The blow-up of \(\mathbb P^2\times\mathbb P^1\) along a fiber conic]
\label{exa:blowup-conic-P2P1}
Let \(X\) be the blow-up of \(\mathbb P^2\times\mathbb P^1\) along \(C\times\{p\}\), where \(C\subset\mathbb P^2\) is a smooth conic and \(p\in\mathbb P^1\).
This is the non-horospherical spherical Fano threefold 3-2-17 of~\cite{Delcroix_Montagard}, under the action of \(G=\mathrm{SL}_2\times\mathbb G_m\).
Writing \(H,F\) for the pullbacks of the hyperplane classes and \(E\) for the exceptional divisor, we consider the anticanonical polarization \(L=-K_X=3H+2F-E\).

Let \(\alpha\) be the simple root of \(\mathrm{SL}_2\) and \(\chi\) a primitive character of \(\mathbb G_m\).
The spherical lattice is \(M=\mathbb Z(2\alpha)\oplus\mathbb Z\chi\) and there is a single spherical root \(\Sigma=\{2\alpha\}\). 
For a choice of linearization, the moment polytope is
\[
\Delta=\left\{
s\alpha+t\chi:
0\leq t\leq2,\quad
0\leq s\leq\min(3,1+2t)
\right\}.
\]
For the complete linear system,
\(
J=\{\alpha,\alpha+\chi,\alpha+2\chi,
3\alpha+\chi,3\alpha+2\chi\}.
\)
Writing \(x=s\alpha+t\chi\), the decomposition polytope is
\[
K_x=\left\{
r\in\mathbb R_{\geq0}^5:
\begin{array}{l}
\sum_{i=0}^4r_i=1,\\
r_1+2r_2+r_3+2r_4=t,\\
r_0+r_1+r_2+3r_3+3r_4\geq s
\end{array}
\right\}. 
\]
The Duistermaat--Heckman measure is \(P_{\mathrm{DH}}\,d\mu_M=s\,ds\,dt\). 
Consequently, \(\Sigma_{\mathrm{sec}}^G(X,L) = 24\int_\Delta K_x\,s\,ds\,dt\).

Eliminating \(r_3\) and \(r_4\), we find that the projection
of \(K_x\) onto its first three coordinates is given by
\[
\begin{cases}
r_0,r_1,r_2\geq0,\\
2r_0+r_1\leq 2-t,\\
r_2-r_0\leq t-1,\\
r_0+r_1+r_2\leq(3-s)/2.
\end{cases}
\]
The chamber complex, pictured below, consists of ten maximal cells,
obtained by cutting \(\Delta\) along the lines
\(t=1\), \(s=1\), \(s=1+t\), and, for \(t\geq1\),
\(s=5-2t\) and \(s=2t-1\).

\begin{figure}
\begin{center}
\begin{tikzpicture}[x=1.35cm,y=.78cm]
\draw[thick]
  (0,0)--(0,1)--(1,3)--(2,3)--(2,0)--cycle;
\draw (1,0)--(1,3);
\draw (0,1)--(2,1);
\draw (0,1)--(2,3);
\draw (1,3)--(2,1);
\draw (1,1)--(2,3);
\draw[->] (0,0)--(2.42,0) node[right] {$t$};
\draw[->] (0,0)--(0,3.45) node[above] {$s$};
\foreach \t in {0,1,2}
  \node[below] at (\t,0) {$\t$};
\foreach \s in {1,2,3}
  \node[left] at (0,\s) {$\s$};
\end{tikzpicture}
\caption{The chamber complex associated with Example~\ref{exa:blowup-conic-P2P1}.}
\label{fig:chamber-complex}
\end{center}
\end{figure}
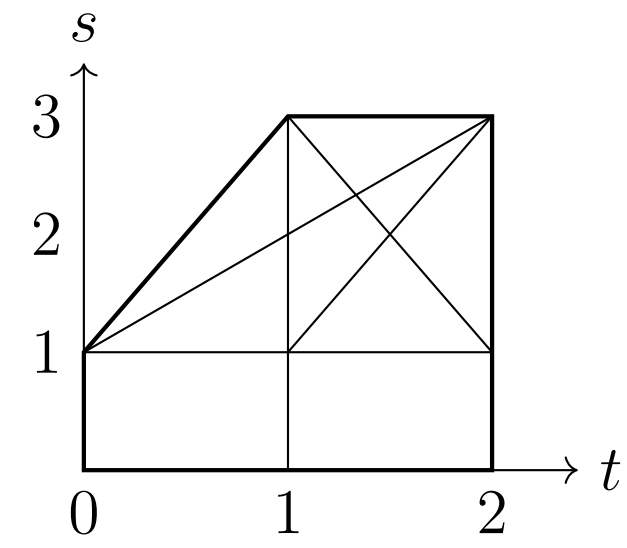

We then compute the spherical secondary polytope using
the Minkowski sum formula of Remark~\ref{rem:BS}.
This is a three-dimensional polytope contained in the affine space defined by 
\[
\sum_{i=0}^4p_i=160,\qquad
p_1+2p_2+p_3+2p_4=196.
\]
Its projection onto the coordinates \((p_0,p_1,p_2)\) is described by
\[
\left\{
(a,b,c)\in\mathbb R^3:
\begin{array}{l}
a\geq18,\quad b\geq0,\quad c\geq0,\\
c-a\leq16,\quad a+b\leq70,\\
2a+b\leq104,\quad a+b+c\leq88
\end{array}
\right\}.
\]
In particular, the spherical secondary polytope has ten vertices, fifteen edges and seven facets as shown in Figure~\ref{fig:secondary-blowup-conic}.
\begin{figure}[ht]
\centering
\begin{tikzpicture}[
  x={(0.100cm,-0.043cm)},
  y={(0.067cm, 0.064cm)},
  z={(0cm,0.092cm)},
  edge/.style={draw=black,line width=0.9pt},
  hidden/.style={draw=gray!65,densely dashed,line width=0.7pt}
]
\coordinate (v0) at (18, 0, 0);
\coordinate (v1) at (18, 0,34);
\coordinate (v2) at (18,52, 0); 
\coordinate (v3) at (18,36,34);
\coordinate (v4) at (18,52,18);
\coordinate (v5) at (52, 0, 0);
\coordinate (v6) at (36, 0,52);
\coordinate (v7) at (52, 0,36);
\coordinate (v8) at (34,36, 0);
\coordinate (v9) at (34,36,18);

\draw[hidden] (v0)--(v2) (v2)--(v4) (v2)--(v8);

\fill[gray!28,fill opacity=.40] (v0)--(v5)--(v7)--(v6)--(v1)--cycle;
\fill[gray!65,fill opacity=.40] (v1)--(v6)--(v3)--cycle;
\fill[gray!65,fill opacity=.35] (v5)--(v8)--(v9)--(v7)--cycle;
\fill[gray!65,fill opacity=.40] (v3)--(v4)--(v9)--(v7)--(v6)--cycle;

\draw[edge] (v0)--(v1)--(v3)--(v4)--(v9)--(v8)--(v5)--(v0);
\draw[edge] (v0)--(v5)--(v7)--(v6)--(v1);
\draw[edge] (v3)--(v6)--(v7)--(v9);
\draw[edge] (v5)--(v7);

\foreach \v in {v0,v1,v3,v4,v5,v6,v7,v8,v9}
  \fill[black] (\v) circle[radius=1.4pt];
\fill[gray!65] (v2) circle[radius=1.1pt];
\end{tikzpicture}
\caption{The spherical secondary polytope of
$\operatorname{Bl}_{C\times \{p\}}(\mathbb P^2\times\mathbb P^1)$
with its anticanonical polarization.}
\label{fig:secondary-blowup-conic}
\end{figure}
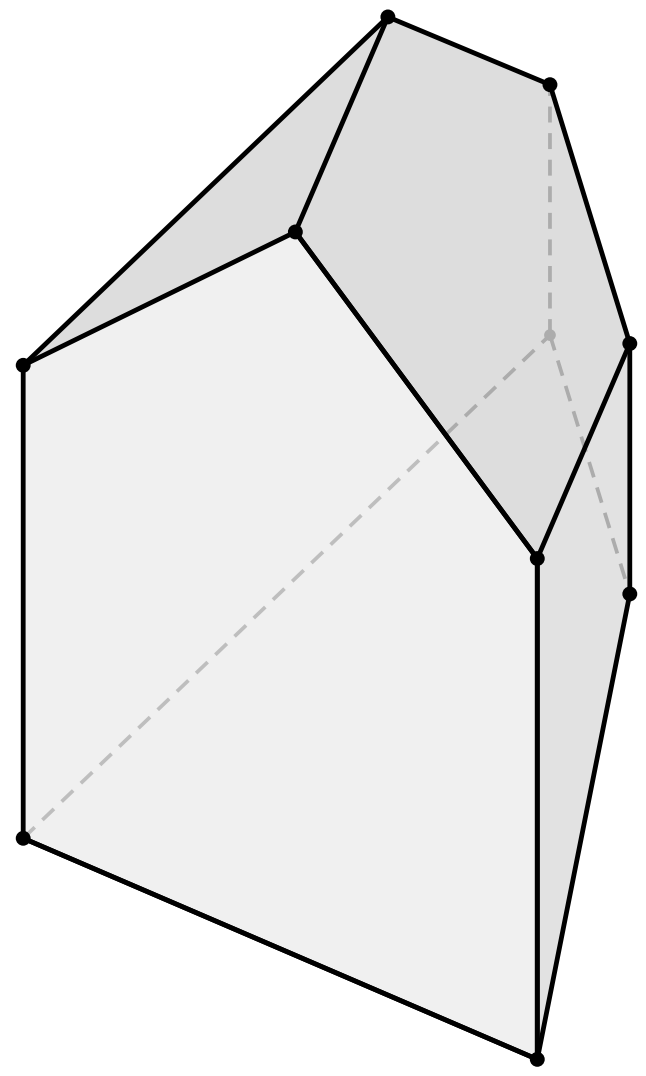
\end{exa}

\section{The equivariant Chow polytope of spherical varieties}
\label{sec:equiv_Chow_polytope}

In this section we identify the equivariant Chow polytope of \(\pi_*[X]\) with the spherical secondary polytope. 
A one-parameter subgroup of the block-scaling torus assigns heights to the irreducible summands of \(E\). 
We show that the associated spherical envelope extends these heights to a filtration of the full section ring and determines the integral closure of the algebra defining the degeneration. 
The Hilbert--Chow weight formula then expresses the weight of the limiting Chow form in terms of the leading coefficient of the total section weight. 
This coefficient is the Duistermaat--Heckman integral of the spherical envelope, 
yielding equality of the two support functions.

\subsection{The Chow-polytope formula}

For \(\lambda\in J\), let \(\varepsilon_\lambda\in X^*(\widetilde T_0)\) be the character of the scalar action on \(W_{1,\lambda}^*\subset E^*\).
These characters identify the block-coordinate space with the real character space through
\[
\iota:\mathbb R^J\longrightarrow
X^*(\widetilde T_0)_{\mathbb R},
\qquad
(r_\lambda)\longmapsto
\sum_{\lambda\in J}r_\lambda\varepsilon_\lambda.
\]
Thus a cocharacter \(\zeta\) with \(h_\lambda=\langle\varepsilon_\lambda,\zeta\rangle\) acts on \(W_{1,\lambda}\subset E\) with weight \(-h_\lambda\). 
We now state the main result, using the notation introduced in Section~\ref{sec:Chow}. 

\begin{thm}\label{thm:spherical-secondary-Chow}
Let \(E\subset H^0(X,L)\) be a base-point-free \(G\)-invariant linear system defining a finite morphism \( \pi:X\longrightarrow Y\subset\mathbb P(E^*) \).
Then
\[
\operatorname{Chow}_{\mathrm{GL}}^G(\pi_*[X])
=
\iota\bigl(\Sigma_{\mathrm{sec}}^G(X,L;E)\bigr).
\]
Writing \(\rho\) for the composition of \(\iota\) with restriction of characters to \(T_0\), one has
\[
\operatorname{Chow}_{\mathrm{SL}}^G(\pi_*[X])
=
\rho\bigl(\Sigma_{\mathrm{sec}}^G(X,L;E)\bigr).
\]
\end{thm}

The remainder of this section is devoted to the proof of Theorem~\ref{thm:spherical-secondary-Chow}, by proving equality of the support functions. 
We first describe the degeneration associated with a one-parameter subgroup of the block-scaling torus.

Let \(\zeta:\mathbb C^*\to\widetilde T_0\) be a cocharacter with weight \(h_\lambda=\langle \varepsilon_{\lambda},\zeta\rangle\) on \(W_{1,\lambda}^*\subset E^*\). 
Using \(\tau=t^{-1}\), consider the family of embedded varieties
\[
\zeta(\tau)^{-1}\cdot Y\subset\mathbb P(E^*),
\qquad \tau\neq0.
\]
Its closure \(\mathcal{Y}_h\) in \(\mathbb P(E^*)\times\mathbb C\) is a degeneration of \(Y\).

Since both support functions are positively homogeneous in the cocharacter, it suffices to prove their equality after replacing each cocharacter by a sufficiently divisible positive multiple. 
We may therefore make a finite base change \(\tau\mapsto\tau^q\) in the degeneration, which leaves the limiting Chow point unchanged and multiplies the heights and their spherical envelope by \(q\). 
Choose \(q\) clearing the denominators of the affine pieces of the spherical envelope \(g_h\).
Retaining the notation \(\tau\), \(h\), and \(g_h\) for the new parameter, heights, and envelope, we thus assume
\[
k\,g_h(\mu/k)\in\mathbb Z
\qquad(k\geq1,\ \mu\in\Lambda_k).
\]

With \(R\) and \(A_E\) as in~\eqref{eqn:defn_R_A_E}, the homogeneous coordinate algebra of this family is
\[
D_h=
\mathbb C\!\left[
\tau,\,
\tau^{-h_\lambda}W_{1,\lambda}
\;:\;\lambda\in J
\right]
\subset R[\tau,\tau^{-1}],
\]
with \(\tau\) of section degree zero. 
Indeed, its degree-one generators scale the projective coordinates in block \(\lambda\) by \(\tau^{-h_\lambda}\).
Moreover, \(D_h[\tau^{-1}]=A_E[\tau,\tau^{-1}]\), accounting for the triviality over \(\mathbb C^*\)

To extend the finite morphism \(X\to Y\) over this degeneration, we take the integral closure of \(D_h\) in \(\operatorname{Frac}(R)(\tau)\).

\begin{lem}\label{lem:finite-map-normalization-rees}
The integral closure of \(D_h\) in \(\operatorname{Frac}(R)(\tau)\) is
\[
\mathcal R_h
=
\mathbb C[\tau]\oplus
\bigoplus_{k\geq1}
\bigoplus_{\mu\in\Lambda_k}
\tau^{-k g_h(\mu/k)}
W_{k,\mu}\,\mathbb C[\tau].
\]
In particular, \(\mathcal R_h\) is finite over \(D_h\), and \(\mathcal R_h[\tau^{-1}]=R[\tau,\tau^{-1}]\).
\end{lem}

\begin{proof}
Let \(\overline D_h\) denote the integral closure of \(D_h\) in the field \(K:=\operatorname{Frac}(R)(\tau)\).
As established in Proposition~\ref{prop:decomposition-polytopes-nonempty}, \(R\) is finite over \(A_E\).
Since \( D_h[\tau^{-1}]=A_E[\tau,\tau^{-1}]\), \(\operatorname{Frac}(D_h)=\operatorname{Frac}(A_E)(\tau)\) and the extension
\(
K/\operatorname{Frac}(D_h)
=
\operatorname{Frac}(R)(\tau)/
\operatorname{Frac}(A_E)(\tau)
\)
is finite.
The algebra \(D_h\) is finitely generated over \(\mathbb C\), so its integral closure in this finite field extension is a finite \(D_h\)-module~\cite[Corollary~3.13, p.293]{Eisenbud_1995}.

The full section ring \(R\) is normal because \(X\) is normal. 
Hence \(R[\tau,\tau^{-1}]\) is normal.
It is also finite over \(A_E[\tau,\tau^{-1}]\), with fraction field \(K\), and therefore is the integral closure of \(A_E[\tau,\tau^{-1}]\) in \(K\).
Since integral closure commutes with localization,
\(
\overline D_h[\tau^{-1}]
=
R[\tau,\tau^{-1}].
\)

Write the spherical envelope as
\[
g_h(x)=\min_{1\leq i\leq m}\ell_i(x),
\qquad
\ell_i(x)=b_i+\langle x-\chi,w_i\rangle,
\quad w_i\in\mathcal V_{\mathbb Q}.
\]

The normal affine variety \(\operatorname{Spec}R[\tau,\tau^{-1}]\) is spherical for \(\mathbb G_m\times G\times\mathbb G_m\), where the two \(\mathbb G_m\) actions record the section degree and the parameter exponent.
Its coordinate ring decomposes into irreducible submodules \(\tau^aW_{k,\mu}\).
On its weight lattice, consider the linear form
\[
L_i(k,\mu,a)
=a+kb_i+\langle\mu-k\chi,w_i\rangle.
\]
By~\eqref{eq:section-ring-multiplication}, every multiplication tail is of the form \((0,\gamma,0)\), with \(\gamma\in\mathcal V^\vee\).
Since \(w_i\in\mathcal V\), we have \(L_i(0,\gamma,0)=\langle\gamma,w_i\rangle\leq0\).
By the characterization of the valuation cone in terms of multiplication tails~\cite[Lemma~6.1]{Knop_1991}, \(L_i\) thus determines an invariant valuation \(v_i\) on \(K\).

The filtration associated with \(v_i\) is invariant under \(\mathbb G_m\times G\times\mathbb G_m\).
By multiplicity-freeness, each filtration subspace is a direct sum of irreducible submodules.
Consequently,
\[
v_i\left(\sum_{k,\mu,a}\tau^a u_{k,\mu,a}\right)
=
\min_{u_{k,\mu,a}\neq0}
\left\{a+kb_i+\langle\mu-k\chi,w_i\rangle\right\},
\]
where \(u_{k,\mu,a}\in W_{k,\mu}\).

Let \(\mathcal O_{v_i}\subset K\) be the corresponding valuation ring. For an integer \(a\),
\[
\begin{aligned}
\tau^aW_{k,\mu}\subset\mathcal O_{v_i}
\quad\text{for every }i
&\iff a+k\ell_i(\mu/k)\geq0
       \quad\text{for every }i\\
&\iff a+k g_h(\mu/k)\geq0.
\end{aligned}
\]
Consequently, by definition, 
\(
\mathcal R_h
=
R[\tau,\tau^{-1}]
\cap\bigcap_{i=1}^m\mathcal O_{v_i}.
\)
Each ring in this intersection is integrally closed in \(K\), so \(\mathcal R_h\) is integrally closed.

We have 
\(
g_h(\lambda)\geq h_\lambda,
\)
for \(\lambda\in J\) and therefore \(D_h\subset\mathcal R_h\). 
Since \(\mathcal R_h\) is integrally closed in \(K\), we obtain
\(
\overline D_h\subset\mathcal R_h.
\)

The finite inclusion \(D_h\subset\overline D_h\) induces a finite morphism
\[
\overline\Pi:\mathcal X
:=\operatorname{Proj}_{\mathbb C[\tau]}\overline D_h
\longrightarrow\mathcal Y_h.
\]
Let \(\mathcal L=\overline\Pi^*\mathcal O_{\mathcal Y_h}(1)\).
The pair \((\mathcal X,\mathcal L)\) is a normal \(G\)-equivariant test configuration for \((X,L)\).

By~\cite[Theorem~4.1 and Sections~4.2--4.3] {Delcroix_2023_KSSV2}, its section filtration is described by a rational concave piecewise-affine function \(f\) on \(\Delta\), with slopes in \(\mathcal V\).
We normalize \(f\) using the trivialization over \(\mathbb C^*\) inherited from \(\overline D_h[\tau^{-1}]=R[\tau,\tau^{-1}]\).

Since the weighted generators \(\tau^{-h_\lambda}W_{1,\lambda}\) belong to \(\overline D_h\), they extend as sections on \(\mathcal X\). 
Hence \(f(\lambda)\geq h_\lambda\) for \(\lambda\in J\).
Every affine piece defining \(f\) is therefore an admissible majorant in Lemma~\ref{lem:spherical-envelope-duality}.
Consequently,
\[
f\geq g_h.
\]

Conversely, since \(\overline D_h\) is a finite graded \(D_h\)-module and \(D_h\) is generated in degree one over \(\mathbb C[\tau]\), the canonical maps
\(
(\overline D_h)_k
\longrightarrow H^0(\mathcal X,\mathcal L^{\otimes k})
\)
are isomorphisms for all sufficiently large \(k\), see ~\cite[Exercise II.5.9(b)]{Hartshorne_1977_book},
applied to the module \(\overline D_h\) on \(\operatorname{Proj}D_h\).
These identifications respect the trivializations over \(\mathbb C^*\).
For each rational \(x\in\Delta\), comparison of filtration levels in sufficiently large divisible degrees, using \(\overline D_h\subset\mathcal R_h\), therefore gives
\[
f(x)\leq g_h(x).
\]
Continuity extends this inequality to all of \(\Delta\). 
Together with the opposite inequality proved above, this yields \(f=g_h\).

The section description then shows that \(\overline D_h\) and \(\mathcal R_h\) agree in all
sufficiently large divisible degrees.
Let \(u\in\mathcal R_h\) be homogeneous of positive section degree. 
For a sufficiently large divisible integer \(m\), we have \(u^m\in\overline D_h\).
Thus \(u\) is integral over \(\overline D_h\), and belongs to it because \(\overline D_h\) is integrally closed in \(K\).
The degree-zero part \(\mathbb C[\tau]\) already lies in \(D_h\).
Therefore \(\mathcal R_h=\overline D_h\), as required.
\end{proof}

Geometrically, the finite inclusion \(D_h\subset\mathcal R_h\) induces a finite morphism
\[
\Pi:\mathcal X_h
:=\operatorname{Proj}_{\mathbb C[\tau]}\mathcal R_h
\longrightarrow
\mathcal Y_h
:=\operatorname{Proj}_{\mathbb C[\tau]}D_h
\]
where the total space \(\mathcal X_h\) is normal. 
Over \(\mathbb C^*\) this morphism identifies abstractly with
\[
\pi\times\mathrm{id}_{\mathbb C^*}:
X\times\mathbb C^*\longrightarrow Y\times\mathbb C^*.
\]
Thus \(\Pi\) extends the original finite morphism over the degeneration of its image. 
If \(E=H^0(X,L)\) and \(L\) is very ample, then \(\Pi\) is simply the normalization morphism for \(\mathcal Y_h\). 

The explicit formula shows that the spherical envelope determines the filtration on the source section ring: the summand \(W_{k,\mu}\) has filtration level \(k\,g_h(\mu/k)\). 
In particular, the special-fiber algebra is 
\[
\mathcal R_h/(\tau)\simeq\operatorname{gr}_{g_h}R.
\]

We will now use the relation between Hilbert weights and Chow weights
established in~\cite[Theorem~3.3]{Kapranov_Sturmfels_Zelevinsky_1992}.
For our purposes, we need its extension from structure sheaves of
projective subschemes to arbitrary equivariant coherent sheaves.

Let \(\mathcal{F}\) be a coherent sheaf on \(\mathbb{P}(V)\). 
Let \(Y_1,\ldots,Y_s\) be the \(n\)-dimensional irreducible components of \(\operatorname{Supp}\mathcal F\). 
For each \(i\), let \(m_i\) be the length of \(\mathcal F_{\eta_i}\) over \(\mathcal O_{\mathbb P(V),\eta_i}\), where \(\eta_i\) is the generic point of \(Y_i\).
The associated cycle is
\[
Z=[\mathcal F]_n=\sum_{i=1}^s m_i[Y_i].
\]
The finite morphism \(\pi\) allows us to represent the cycle \(\pi_*[X]\) by the coherent sheaf \(\pi_*\mathcal O_X\) on the image.

\begin{lem}\label{lem:coherent-Hilbert-Chow}
Let \(\zeta:\mathbb C^*\to\mathrm{GL}(V)\), let \(\mathcal F\) be a \(\zeta\)-equivariant coherent sheaf on \(\mathbb P(V)\) whose support has dimension \(n\), and let \(Z=[\mathcal F]_n\) denote the associated cycle. 
For \(k\gg0\), write the total weight as
\[
w_{\mathcal F}(k)
:=
\operatorname{wt}_{\zeta}
\det H^0(\mathbb P(V),\mathcal F(k))
=
b_{n+1}k^{n+1}+O(k^n).
\]
Then the weight of the Chow form of \(Z\), under the canonical action on the Chow representation, is 
\[
\operatorname{wt}_{\zeta}(R_Z)
=
-(n+1)!\,b_{n+1}.
\]
\end{lem}

\begin{proof}
First assume that \(\mathcal F=\mathcal O_Y\) for a \(\zeta\)-invariant irreducible subvariety \(Y\subseteq\mathbb P(V)\) of dimension \(n\).
By~\cite[Theorem~3.3]{Kapranov_Sturmfels_Zelevinsky_1992}, applied with respect to a maximal torus of \(\mathrm{GL}(V)\) containing the image of \(\zeta\), we have
\[
\operatorname{wt}_\zeta(R_Y)
=
-(n+1)!\lim_{k\to\infty}
\frac{w_{\mathcal O_Y}(k)}{k^{n+1}}.
\]
The minus sign comes from the convention that the Hilbert point is represented by the dual of the determinant of the degree-\(k\) quotient of the homogeneous coordinate ring of \(\mathbb P(V)\).

For the general case, choose a finitely generated equivariant graded module over \(\operatorname{Sym}(V^*)\) representing \(\mathcal F\).
An equivariant homogeneous prime filtration of this module induces a filtration of \(\mathcal F\) whose successive quotients are structure sheaves of invariant irreducible subvarieties, up to twists and characters of \(\mathbb C^*\).
For each \(n\)-dimensional component \(Y_i\) of \(\operatorname{Supp}\mathcal F\), the number of such quotients supported on \(Y_i\) is precisely the generic length \(m_i\).

By Serre vanishing, total weights are additive in equivariant short exact sequences for sufficiently large \(k\).
Furthermore, a coherent sheaf supported in dimension at most \(n-1\) contributes \(O(k^n)\) to the total weight, since its space of sections has dimension \(O(k^{n-1})\) and its weights grow at most linearly in \(k\).
Neither a fixed twist nor a change of linearization by a character affects the coefficient of \(k^{n+1}\).
We thus obtain
\[
w_{\mathcal F}(k)
=
\sum_{i=1}^s m_i w_{\mathcal O_{Y_i}}(k)
+
O(k^n).
\]
Since
\(
R_Z=\prod_{i=1}^s R_{Y_i}^{m_i},
\)
the desired equality follows from the case of structure sheaves.
\end{proof}

We now return to the finite morphism of families \( \Pi:\mathcal X_h\longrightarrow\mathcal Y_h \subset\mathbb P(E^*)\times\mathbb C \).
Set \( \mathcal F=\Pi_*\mathcal O_{\mathcal X_h} \), regarded as a coherent sheaf on \(\mathbb P(E^*)\times\mathbb C\) via the closed embedding.
Let \(\mathcal F_0\) denote its restriction to \(\mathbb P(E^*)\times\{0\}\), identified with \(\mathbb P(E^*)\).

\begin{lem}\label{lem:finite-pushforward-Chow-limit}
The sheaf \(\mathcal F\) is flat over \(\mathbb C\).
Moreover, its special fiber satisfies
\[
[\mathcal F_0]_n
=
\lim_{t\to\infty}\zeta(t)\cdot\pi_*[X],
\]
where the limit is taken in the Chow variety of \(n\)-cycles in \(\mathbb P(E^*)\).
\end{lem}

\begin{proof}
The variety \(\mathcal X_h\) is integral and dominates \(\mathbb C\). 
Hence multiplication by any nonzero polynomial in \(\tau\) is injective on \(\mathcal O_{\mathcal X_h}\), and therefore on its finite pushforward \(\mathcal F\).
Thus \(\mathcal F\) is torsion-free over \(\mathbb C[\tau]\). 
Since \(\mathbb C[\tau]\) is a principal ideal domain, \(\mathcal F\) is flat over \(\mathbb C\).

Over \(\mathbb C^*\), the morphism \(\Pi\) identifies abstractly with \(\pi\times\mathrm{id}_{\mathbb C^*}\), while the image is embedded in the fixed projective space by \(\zeta(\tau)^{-1}\).
Consequently, for \(\tau\neq0\),
\[
[\mathcal F_\tau]_n
=
\zeta(\tau)^{-1}\cdot\pi_*[X].
\]
Here we use \( [\pi_*\mathcal O_X]_n=\pi_*[X]\): both cycles equal \(e[Y]\), where \(e=[\mathbb C(X):\mathbb C(Y)]\).

Since \(\Pi\) is finite, its pushforward commutes with restriction to the special fiber, and
\[
[\mathcal F_0]_n=(\Pi_0)_*[\mathcal X_{h,0}].
\]
The special fibers are the Cartier divisors cut out by \(\tau\).
The projection formula for Cartier divisors~\cite[Proposition~2.3(c)]{Fulton_1998} therefore gives
\[
(\Pi_0)_*[\mathcal X_{h,0}]=e[\mathcal Y_{h,0}].
\]
The right-hand side is the specialization of the generic-fiber cycle, since \(\mathcal Y_h\) is the closure of its generic fiber.
In the fixed projective space, this specialization is the limit in the Chow variety. 
Therefore
\[
[\mathcal F_0]_n
=
\lim_{\tau\to0}\zeta(\tau)^{-1}\cdot\pi_*[X]
=
\lim_{t\to\infty}\zeta(t)\cdot\pi_*[X].
\]
\end{proof}

In the moduli framework of Alexeev and Brion, this specialization is the value on the special fiber of the cycle map~\cite[Proposition~4.2.2]{Alexeev_Brion_2006}, which sends a stable spherical variety finite over projective space to its pushforward cycle.

\begin{proof}[Proof of Theorem~\ref{thm:spherical-secondary-Chow}]
We compare support functions.
Fix a cocharacter \(\zeta:\mathbb C^*\to\widetilde T_0\), and put \( h_\lambda=\langle\varepsilon_\lambda,\zeta\rangle\).
By positive homogeneity, we may replace \(\zeta\) by a sufficiently divisible positive multiple and assume that
\[
p_{k,\mu}:=k\,g_h(\mu/k)\in\mathbb Z
\qquad(k\geq1,\ \mu\in\Lambda_k).
\]

Use the families \(\mathcal X_h\), \(\mathcal Y_h\) and the sheaf \(\mathcal F\) constructed above.
We equip the image family with the action \(a\cdot(y,\tau)=(\zeta(a)y,a^{-1}\tau)\).
This action lifts to \(\mathcal X_h\).
Under the trivialization over \(\mathbb C^*\), it acts trivially on \(X\) and by inverse scaling on the parameter.
Consequently, its action on \(R[\tau,\tau^{-1}]\) satisfies \(a\cdot u=u\) for \(u\in R\) and \(a\cdot\tau=a\tau\). 
In particular, for \(u\in W_{k,\mu}\), 
\[
a\cdot(\tau^{-p_{k,\mu}}u)
=
a^{-p_{k,\mu}}\tau^{-p_{k,\mu}}u.
\]
Thus the summand \(W_{k,\mu}\) surviving in \((\mathcal R_h/(\tau))_k\) has weight \(-p_{k,\mu}\).

Since \(\Pi\) is finite, passage to the fiber over \(\tau=0\) gives \(\mathcal F_0 = (\Pi_0)_*\mathcal O_{\mathcal X_{h,0}}\). 
As a sheaf on \(\mathcal Y_{h,0}\), it is the sheaf associated with the finite graded
\(D_h/(\tau)\)-module \(\mathcal R_h/(\tau)\). 
Consequently, for \(k\gg0\), we have an equivariant identification \(H^0\left(\mathbb P(E^*),\mathcal F_0(k)\right) \simeq (\mathcal R_h/(\tau))_k\).
The total weight is therefore
\[
w_{\mathcal F_0}(k)
=
-\sum_{\mu\in\Lambda_k}
k\,g_h(\mu/k)\dim W_{k,\mu}.
\]

The Weyl dimension formula, followed by lattice Riemann-sum asymptotics on the finitely many polyhedral cells on which \(g_h\) is affine (see~\cite[Section~1.4]{Pukhlikov_Khovanskii_1992}), gives
\[
\begin{aligned}
w_{\mathcal F_0}(k)
&=
-k^{n-r+1}
\sum_{\mu\in\Lambda_k}
g_h(\mu/k)P_{\mathrm{DH}}(\mu/k)
+O(k^n)\\
&=
-k^{n+1}\int_\Delta
g_h(x)P_{\mathrm{DH}}(x)\,d\mu_M(x)
+O(k^n).
\end{aligned}
\]

Set \(Z=\pi_*[X]\) and \(Z_0=[\mathcal F_0]_n\).
By Lemma~\ref{lem:finite-pushforward-Chow-limit}, \(Z_0=\lim_{t\to\infty}\zeta(t)\cdot Z\).
Thus the Chow line of \(Z_0\) is the maximal-weight limit of the Chow line of \(Z\), and
\[
\max_{\xi\in\operatorname{Chow}_{\mathrm{GL}}^G(Z)}
\langle\xi,\zeta\rangle
=
\operatorname{wt}_\zeta(R_{Z_0}).
\]
Applying Lemma~\ref{lem:coherent-Hilbert-Chow} to \(\mathcal F_0\), we obtain
\[
\operatorname{wt}_\zeta(R_{Z_0})
=(n+1)!\int_\Delta
g_h(x)P_{\mathrm{DH}}(x)\,d\mu_M(x).
\]
By the definition of \(g_h\) and of the Minkowski integral, this is precisely
\(s_{\iota(\Sigma_{\mathrm{sec}}^G(X,L;E))}(\zeta)\).
Positive homogeneity removes the divisibility assumption on \(\zeta\).
Equality on integral cocharacters extends to rational cocharacters by homogeneity, and then to all real cocharacters by continuity.
The two polytopes therefore have the same support function, proving the first assertion. 
Restricting characters to \(T_0\) gives the second assertion.
\end{proof}

\section{Chow stability of spherical varieties}
\label{sec:CSSV}

We now apply Theorem~\ref{thm:spherical-secondary-Chow} to the study of Chow stability of polarized spherical varieties.
We first derive a criterion for Chow semistability and polystability in terms of convex functions on the moment polytope.
We then investigate the asymptotic behavior of the corresponding Chow functional, making explicit its relation to the K-stability criterion and deriving necessary conditions for asymptotic Chow semistability.

\subsection{Chow stability in terms of convex functions}

Let \((X,L)\) be a polarized spherical variety of dimension \(n\).
Fix \(k>0\) such that \(L^{\otimes k}\) is very ample, and consider the complete embedding \(X\hookrightarrow\mathbb P(H^0(X,L^{\otimes k})^*)\).
Write \(H^0(X,L^{\otimes k}) = \bigoplus_{\lambda\in\Lambda_k}W_{k,\lambda}\), 
and set 
\[ D_k=\deg(X,L^{\otimes k}) \qquad p_k=\sum_{\lambda\in\Lambda_k}\dim W_{k,\lambda}.\]
Define the vector
\[
\mathbf b_k
=
\frac{(n+1)D_k}{p_k}
\bigl(\dim W_{k,\lambda}\bigr)_{\lambda\in\Lambda_k}
\in\mathbb R^{\Lambda_k}.
\]

\begin{prop}\label{prop:spherical-Chow-stability}
For the complete embedding defined by \(L^{\otimes k}\),
\[
\begin{aligned}
X\text{ is Chow semistable}
&\iff
\mathbf b_k\in
\Sigma_{\mathrm{sec}}^G(X,L^{\otimes k}),\\
X\text{ is Chow polystable}
&\iff
\mathbf b_k\in
\operatorname{relint}
\Sigma_{\mathrm{sec}}^G(X,L^{\otimes k}).
\end{aligned}
\]
\end{prop}

\begin{proof}
Let \(\rho_k:\mathbb R^{\Lambda_k}\longrightarrow
X^*(T_{0,k})_{\mathbb R}\) be the restriction of block
characters to the connected centralizer \(T_{0,k}\)
of \(G\) in \(\mathrm{SL}(H^0(X,L^{\otimes k})^*)\).
The determinant character gives
\[
\ker\rho_k
=
\mathbb R\bigl(\dim W_{k,\lambda}\bigr)_{\lambda\in\Lambda_k}.
\]
Moreover, since every decomposition vector has coordinate sum \(1\), the definition of the spherical secondary polytope and the degree formula imply
\[
\Sigma_{\mathrm{sec}}^G(X,L^{\otimes k})
\subset
\left\{u:\sum_{\lambda\in\Lambda_k}u_\lambda
=(n+1)D_k\right\}.
\]
This hyperplane meets \(\ker\rho_k\) at the unique point \(\mathbf b_k\).

By Theorem~\ref{thm:spherical-secondary-Chow},
\(
\operatorname{Chow}_{\mathrm{SL}}^G(X)
=
\rho_k\bigl(\Sigma_{\mathrm{sec}}^G(X,L^{\otimes k})\bigr).
\)
Hence \(0\in\operatorname{Chow}_{\mathrm{SL}}^G(X)\) if and only if \(\mathbf b_k\) belongs to the spherical secondary polytope.
Likewise, since linear maps preserve relative interiors, \(0\in\operatorname{relint}\operatorname{Chow}_{\mathrm{SL}}^G(X)\) if and only if \(\mathbf b_k\) belongs to its relative interior.
The assertions follow from the equivariant Chow semistability and polystability criteria.
\end{proof}

\begin{exa}
For the complete embedding of the quadric in Example~\ref{exa:two-block-quadric}, let \(p=\dim U\) and \(q=\dim W\).
We have \(n=p+q-2\), \(D_1=2\), and the two irreducible blocks have dimensions \(p\) and \(q\).
Thus
\[
\mathbf b_1
=\frac{2(p+q-1)}{p+q}(p,q).
\]
On the other hand, we computed \(\Sigma_{\mathrm{sec}}^G(X,L) = \operatorname{conv}\{(2p-2,2q),(2p,2q-2)\}\).
Since
\[
\mathbf b_1
=\frac{p}{p+q}(2p-2,2q)
+\frac{q}{p+q}(2p,2q-2),
\]
it belongs to the relative interior of the spherical secondary polytope. 
Proposition~\ref{prop:spherical-Chow-stability} thereby recovering the Chow polystability of the quadric in its standard embedding. 
\end{exa}

\begin{exa}
Consider the anticanonical embedding of the blow-up in Example~\ref{exa:blowup-diagonal}.
Let \(\varepsilon_1,\varepsilon_2,\varepsilon_3\) be the characters of \(\widetilde T_0\) corresponding to the three irreducible summands, of dimensions \(35,27,10\), respectively.

Since \(n=4\), \(D_1=336\) and \(p_1=72\), the distinguished character is
\[
\mathbf b_1
=\frac{70}{3}
(35\varepsilon_1+27\varepsilon_2+10\varepsilon_3).
\]
By Example~\ref{exa:blowup-diagonal}, the equivariant Chow polytope is the convex hull of
\[
\begin{aligned}
\nu_1&=1680\varepsilon_1,\\
\nu_2&=606\varepsilon_1+1074\varepsilon_2,\\
\nu_3&=606\varepsilon_1+900\varepsilon_2
       +174\varepsilon_3,\\
\nu_4&=1056\varepsilon_1+624\varepsilon_3.
\end{aligned}
\]
One easily checks that \(\mathbf b_1\) belongs to the interior of the triangle \(\operatorname{conv}\{\nu_1,\nu_2,\nu_4\}\), hence to the relative interior of
\(\operatorname{Chow}_{\mathrm{GL}}^G(X)\).
Proposition~\ref{prop:spherical-Chow-stability}
therefore implies Chow polystability.
\end{exa}

We may formulate the criterion in terms of rational convex piecewise-affine functions, in a form similar to the K-polystability criterion~\cite{Delcroix_2023_KSSV2}. 
For a function \(f:\Delta\to\mathbb R\), set
\[
\mathcal C_k(f)
=
\sum_{\lambda\in\Lambda_k}
f(\lambda/k)\frac{\dim W_{k,\lambda}}{p_k}
-
\int_\Delta
f(x) \frac{P_{\mathrm{DH}}(x)\, d\mu_M(x)}
{\int_\Delta
P_{\mathrm{DH}}(x)\,d\mu_M(x)}.
\]

\begin{thm}\label{thm:spherical-Chow-convex-criterion}
The variety \(X\) is Chow semistable if and only if \(\mathcal C_k(f)\geq0\) 
for every rational convex piecewise-affine function \(f\) on \(\Delta\) whose affine pieces have slopes in \(-\mathcal V\).
The variety \(X\) is Chow polystable if and only if these inequalities hold and, among the same functions, \(\mathcal C_k(f)=0\) if and only if \(f\) is affine with slope in \(\operatorname{Lin}(\mathcal V)\).
\end{thm}

\begin{rem}
It suffices to test the criterion on negatives of spherical envelopes.
More precisely, for \(h\in\mathbb R^{\Lambda_k}\), let \(g_h^{(k)}\) be the spherical envelope on \(k\Delta\) associated with the degree-\(k\) configuration \(\Lambda_k\).
One may restrict the test functions to the functions 
\[
f_h(x)=-g_h^{(k)}(kx),
\qquad x\in\Delta.
\]
These functions are convex and piecewise-affine,
with slopes in \(-\mathcal V\).
\end{rem}

\begin{rem}
We can alternatively consider all continuous convex functions \(f:\Delta\to\mathbb R\) which are supremum of affine functions with slopes in \(-\mathcal V\) as test functions. 
\end{rem}

\begin{proof}
For a function \(u\) on \(\Delta\), write \(\mathcal C_k(u)=A_k(u)-A(u)\) where 
\[
A(u)=
\int_\Delta u\frac{P_{\mathrm{DH}}\,d\mu_M}
     {\int_\Delta P_{\mathrm{DH}}\,d\mu_M},
\qquad
A_k(u)=
\sum_{\lambda\in\Lambda_k}
u(\lambda/k)\frac{\dim W_{k,\lambda}}{p_k}.
\]
For a height vector \(h\), also write
\[
B_k(h)=
\sum_{\lambda\in\Lambda_k}
h_\lambda \frac{\dim W_{k,\lambda}}{p_k},
\qquad
g(x)=g_h^{(k)}(kx).
\]
By the support-function formula for the spherical secondary polytope, the change of variables \(y=kx\), and Brion's degree formula, we have
\[
s_{\Sigma_{\mathrm{sec}}^G(X,L^{\otimes k})}(h)
=(n+1)D_k A(g).
\]
On the other hand, by definition,
\[
\langle\mathbf b_k,h\rangle=(n+1)D_k B_k(h).
\]
Proposition~\ref{prop:spherical-Chow-stability} thus shows that Chow semistability is equivalent to
\[
H_k(h):=A(g)-B_k(h)\geq0
\]
for every rational height vector \(h\), since membership in a compact convex set is characterized by the corresponding support-function inequalities.

Suppose first that \(X\) is Chow semistable.
Let \(f\) be an admissible rational convex piecewise-affine function, and set \(h_\lambda=-f(\lambda/k)\).
By the affine-majorant description of the spherical envelope, \(g\leq-f\). 
Hence
\[
0\leq H_k(h)
=A(g)+A_k(f)
\leq A_k(f)-A(f)
=\mathcal C_k(f).
\]
Conversely, suppose that all the stated inequalities hold. 
Given rational heights \(h\), the function \(-g\) is an admissible test function. 
Moreover, \(g(\lambda/k)\geq h_\lambda\), so
\[
H_k(h)
=
\mathcal C_k(-g)
+
\frac1{p_k}\sum_{\lambda\in\Lambda_k}
\dim W_{k,\lambda}
\bigl(g(\lambda/k)-h_\lambda\bigr)
\geq0.
\]
This proves the semistability assertion.

We next record the equality condition for the support-function criterion. 
Put
\[
S_k=\Sigma_{\mathrm{sec}}^G(X,L^{\otimes k}).
\]
A rational height vector \(h\) defines a linear functional constant on \(S_k\) if and only if, after clearing denominators, its block-scaling cocharacter fixes the Chow point, hence preserves the embedded variety \(X\).
Its normalized degeneration is then a product.
By the classification of spherical test configurations~\cite[Theorem~4.1]{Delcroix_2023_KSSV2}, the corresponding envelope \(g\) is affine with slope in \(\operatorname{Lin}(\mathcal V)\).
Since the embedding uses the complete linear system, its degree-one weights satisfy \(g(\lambda/k)=h_\lambda\) for \(\lambda\in\Lambda_k\).
Conversely, heights obtained by restricting such an affine function define graded automorphisms of the section ring: their slopes annihilate every tail. 
They therefore preserve \(X\), and their functionals are constant on \(S_k\).

Assume now that \(X\) is Chow polystable and that \(\mathcal C_k(f)=0\).
For the heights \(h_\lambda=-f(\lambda/k)\), the inequalities above force \(H_k(h)=0\), and \(A(-f-g)=0\).
The function \(-f-g\) is continuous and nonnegative, and the Duistermaat--Heckman measure has full support on \(\Delta\). 
Thus \(g=-f\).
Since \(\mathbf b_k\in\operatorname{relint}S_k\), equality in the support-function inequality implies that \(h\) is constant on \(S_k\).
The previous characterization shows that \(f\) is affine with slope in \(\operatorname{Lin}(\mathcal V)\).

Conversely, suppose the inequalities and their stated equality characterization hold.
If \(H_k(h)=0\), both nonnegative terms in its displayed decomposition vanish. 
Consequently, \(\mathcal C_k(-g)=0\), and \(g(\lambda/k)=h_\lambda\) for \(\lambda\in\Lambda_k\).
The equality hypothesis makes \(g\) affine with slope in \(\operatorname{Lin}(\mathcal{V})\). 
Hence \(h\) is constant on \(S_k\).
The relative-interior criterion now gives Chow polystability. 
It suffices here to consider rational heights, since \(S_k\) and \(\mathbf b_k\) are rational.

Finally, under semistability, every affine function with slope in \(\operatorname{Lin}(\mathcal{V})\) has equality: both \(f\) and \(-f\) are admissible, so their inequalities force \(\mathcal C_k(f)=0\).
This completes the proof of the theorem.

The same argument proves the two remarks.
For a continuous convex function expressed as a supremum of affine functions with slopes in
\(-\mathcal V\), each negative affine function is an admissible majorant of the sampled heights.
Thus \(g\leq-f\) still holds, and the argument for semistability applies without change.

For the equality case with real heights, the subspace of functionals constant on \(S_k\) is rational.
Its characterization by restrictions of affine functions with slopes in \(\operatorname{Lin}(\mathcal{V})\) therefore extends from rational to real heights.
The previous equality argument consequently applies to continuous admissible functions as well.

Conversely, in both the semistability and polystability arguments, it suffices to test the functions \(-g_h^{(k)}(k\,\cdot)\), as the decomposition of \(H_k(h)\) shows.
\end{proof}

\subsection{Asymptotic Chow stability}

\begin{defn}
A polarized variety \((X,L)\) is \emph{asymptotically Chow semistable} (respectively, \emph{asymptotically Chow polystable}) if, for every sufficiently large integer \(k\), the embedding defined by the complete linear system \(H^0(X,L^{\otimes k})\) is Chow semistable (respectively, Chow polystable).
\end{defn}

Let \(Q\) be the coefficient of \(k^{n-r-1}\) in the Weyl dimension polynomial:
\[
\dim W_{k,\lambda}
=
k^{n-r}P_{DH}(\lambda/k)
+k^{n-r-1}Q(\lambda/k)
+O(k^{n-r-2}),
\]
with \(Q=0\) if \(n=r\). 
Set  now 
\[
A(f)=\int_\Delta fP_{DH}\,d\mu_M,
\qquad
B(f)=\frac12\int_{\partial\Delta}fP_{DH}\,d\sigma_M
+\int_\Delta fQ\,d\mu_M,
\]
where \(d\sigma_M\) is the Lebesgue measure on each facet normalized by its tangent lattice.

\begin{prop}
\label{prop:Chow-weight-expansion}
Let \(f\) be a continuous rational piecewise-affine function on \(\Delta\). 
Choose a positive integer \(q\) such that \(q(u-\chi)\in M\) for every vertex \(u\) of a rational polyhedral subdivision on which \(f\) is affine. 
As \(k\to\infty\) through multiples of \(q\),
\[
\mathcal C_k(f)
=
\frac{1}{kA(1)}
\left(
B(f)-\frac{B(1)}{A(1)}A(f)
\right)
+O_f(k^{-2}).
\]
\end{prop}

The latter notation \(O_f(k^{-2})\) means that the control on the remainder depends on \(f\) in general. 

\begin{proof}
Apply the weighted Euler--Maclaurin formula on each cell of the subdivision~\cite{Pukhlikov_Khovanskii_1992,Berline_Vergne_2007}. 
Correcting for repeated counting along common faces, the internal facet terms cancel because \(f\) is continuous. 
Combining the resulting expansion with the Weyl dimension formula gives
\[
\sum_{\lambda\in\Lambda_k}
(\dim W_{k,\lambda})f(\lambda/k)
=
k^nA(f)+k^{n-1}B(f)+O_f(k^{n-2}).
\]
Taking \(f=1\) yields the corresponding expansion of \(p_k\). 
Dividing these expansions and subtracting \(A(f)/A(1)\) proves the assertion.
\end{proof}

Set
\[
\mathcal L(f)
=
2B(f)-2\frac{B(1)}{A(1)}A(f).
\]
By~\cite[Theorem~1.1]{Delcroix_2023_KSSV2} \((X,L)\) is \(G\)-equivariantly K-polystable if and only if \( \mathcal L(f)\geq0 \) for every rational convex piecewise-affine function with slopes in \(-\mathcal V\), with equality precisely for affine functions whose slope lies in \(\operatorname{Lin}(\mathcal V)\). 

The previous expansion becomes
\[
\mathcal C_k(f)
=
\frac{\mathcal L(f)}{2kA(1)}
+O_f(k^{-2}).
\]
Thus asymptotic Chow semistability implies \(G\)-equivariant K-semistability. 
This gives a direct proof in the spherical setting of the equivariant implication, which also follows from the general result of Ross and Thomas~\cite{Ross-Thomas_2007}.
This limit does not give the equality condition for K-polystability: when \(\mathcal L(f)=0\), the sign of \(\mathcal C_k(f)\) depends on higher-order terms. 
We will focus on the case of rank-one spherical varieties in the next section to investigate this issue. 

Asymptotic Chow semistability also imposes the vanishing of the Chow weights of \(G\)-equivariant product test configurations at every level.
Indeed, for \(v\in\operatorname{Lin}(\mathcal V)\), set \(f_v(x)=\langle x-\chi,v\rangle\).
For integral \(v\), the quantity \(\mathcal C_k(f_v)\) is the normalized Chow weight of the corresponding product test configuration.
Since both \(f_v\) and \(-f_v\) are admissible, Chow semistability at level \(k\) implies
\(\mathcal C_k(f_v)=0\).
In the toric case, this was a key observation by Ono~\cite{Ono_2011,Ono_2013}.
The following proposition shows that vanishing at \(n\) distinct levels implies vanishing at every positive level.

\begin{prop}
\label{prop:product-Chow-weight-finite-test}
Let \(X\) be a normal projective spherical variety of dimension \(n\geq1\), and let \(L\) be an ample \(G\)-linearized line bundle.
If \(\mathcal C_{k_j}(f_v)=0\) for every \(v\in\operatorname{Lin}(\mathcal V)\) at \(n\) distinct positive integers \(k_1,\ldots,k_n\), where \(f_v(x)=\langle x-\chi,v\rangle\), then \(\mathcal C_k(f_v)=0\) for every \(k\geq1\) and \(v\in\operatorname{Lin}(\mathcal V)\).
\end{prop}

\begin{proof}
Fix \(v\in\operatorname{Lin}(\mathcal V)\). 
By linearity, we may assume that \(v\) is integral. 
Since \(v\) annihilates every multiplication tail, assigning weight \(\langle\lambda-k\chi,v\rangle\) to \(W_{k,\lambda}\) defines a \(\mathbb C^*\)-action on the section ring. 
Let \(w_v(k)\) be its total weight in degree \(k\).

By higher-cohomology vanishing~\cite[Proposition~1.3.1(vii)]{Alexeev_Brion_2006} and equivariant Riemann--Roch, \(p_k\) and \(w_v(k)\) are polynomial for \(k\geq0\), of degrees at most \(n\) and \(n+1\), respectively. 
By definition, for \(k\geq1\),
\[
k p_k\,\mathcal C_k(f_v)
=w_v(k)-k p_k\,\frac{A(f_v)}{A(1)}.
\]
Proposition~\ref{prop:Chow-weight-expansion} gives \(\mathcal C_k(f_v)=O(k^{-1})\) along sufficiently divisible \(k\). 
Hence the polynomial on the right has degree at most \(n\). 
It vanishes at \(k=0\), since \(w_v(0)=0\), and at the \(n\) positive integers in the hypothesis. 
It is therefore identically zero.
\end{proof}

\section{Chow and K-stability of rank-one spherical varieties}
\label{sec:rk1}

\subsection{Notation}

Let \((X,L)\) be a polarized spherical variety of dimension \(n\) and spherical rank equal to one, and let \(\chi\) be the weight of a rational \(B\)-eigensection of \(L\). 
If \(X\) is non-horospherical, let \(\xi\) be its spherical root. 
If \(X\) is horospherical, choose a primitive generator \(\xi\) of \(M\).

Write \(\lambda=\chi+x\xi\) for \(\lambda\in\chi+M_{\mathbb R}\).
In the coordinate \(x\), the moment polytope \(\Delta\) is an interval \([a,b]\), and the lattice-normalized measure \(d\mu_M\) is \(dx\).

For \(k\geq 1\), the highest-weights occurring in \(H^0(X,L^{\otimes k})\) are
\[
\Lambda_k
=
\{k\chi+j\xi : j\in\mathbb Z,\ ka\leq j\leq kb\}.
\]
The normalized weight \((k\chi+j\xi)/k\) has coordinate \(j/k\).
Thus the coordinates of the normalized weights form the set
\[
\Gamma_k=[a,b]\cap k^{-1}\mathbb Z.
\]

If \(L^{\otimes k}\) is very ample, Proposition~\ref{prop:decomposition-polytopes-nonempty} applied to its complete linear system gives \( [a,b]\subseteq\operatorname{conv}(\Gamma_k)-C\), 
where \(C\) denotes the tail cone in the coordinate \(x\).
Since \(C=\mathbb R_{\geq 0}\) in the non-horospherical case and \(C=\{0\}\) in the horospherical case, we have \(b\in\Gamma_k\) in both cases and \(a\in\Gamma_k\) in the horospherical case.

By Theorem~\ref{thm:spherical-Chow-convex-criterion}, the Chow test functions are the rational convex piecewise-affine functions on \([a,b]\), with the additional requirement that they be nondecreasing if \(X\) is non-horospherical. 
The equality cases allowed by Chow polystability are the constant functions in the non-horospherical case and the affine functions in the horospherical case.

Write \(P(x)=P_{\mathrm{DH}}(\chi+x\xi)\), and denote by \(Q(x)\) the coefficient of \(k^{n-2}\) in the Weyl dimension polynomial in this coordinate, with \(Q=0\) if \(n=1\). 
We recall the notation 
\[
A(f)=\int_a^b f(x)P(x)\,dx,
\qquad
B(f)=\frac{f(a)P(a)+f(b)P(b)}{2}
     +\int_a^b f(x)Q(x)\,dx.
\]
Writing \(p_k=\dim H^0(X,L^{\otimes k})\), the Chow functional is
\[
\mathcal C_k(f)
=
\sum_{j/k\in\Gamma_k}
f(j/k)\frac{\dim W_{k,k\chi+j\xi}}{p_k}
-\frac{A(f)}{A(1)}.
\]

For \(t\in[a,b]\), define the hinge function \(h_t\colon[a,b]\to\mathbb R\) by
\[
h_t(x)=\max\{0,x-t\}.
\]
This is a convex, piecewise-affine function increasing with (at most) two linearity domains. 

\subsection{Finite test for Chow stability}

The following shows that Chow stability at a fixed level is checked by a finite number of tests. 
The proof reduces the admissible function to finitely many hinge functions and then uses strict concavity between consecutive grid points.

\begin{prop}\label{prop:rank-one-finite-Chow-tests}
Suppose that \(L^{\otimes k}\) is very ample.
\begin{enumerate}
\item If \(X\) is non-horospherical, the embedding defined by \(L^{\otimes k}\) is Chow semistable if and only if
\[
\mathcal C_k(h_t)\geq 0
\qquad\text{for every }t\in\Gamma_k\setminus\{b\}.
\]
It is Chow polystable if and only if all these inequalities are strict.

\item If \(X\) is horospherical, the embedding defined by \(L^{\otimes k}\) is Chow semistable if and only if 
\[
\mathcal C_k(x)=0
\qquad\text{and}\qquad
\mathcal C_k(h_t)\geq 0
\quad\text{for every }t\in\Gamma_k\cap(a,b).
\]
It is Chow polystable if and only if \(\mathcal C_k(x)=0\) and all these inequalities are strict.
\end{enumerate}
\end{prop}

\begin{proof}
Modulo constants, every admissible function in the non-horospherical case is a nonnegative linear combination of the \(h_t\) with \(a\leq t<b\). 
In the horospherical case, the same holds modulo affine functions using \(a<t<b\), and testing \(x\) and \(-x\) forces \(\mathcal C_k(x)=0\).

The map \(t\mapsto\mathcal C_k(h_t)\) is strictly concave between consecutive grid points, since its discrete term is affine there and its second derivative is \(-P(t)/\int_a^b P(x)\,dx<0\).
Its values at \(a\) and \(b\) are \(\mathcal C_k(x)\) and \(0\), respectively. 
If \(a<\min\Gamma_k\), it is strictly decreasing on \([a,\min\Gamma_k]\), since its derivative there is \(-\int_a^t P(x)\,dx/\int_a^b P(x)\,dx\) in the interior.
Thus the grid inequalities imply the corresponding inequalities for every hinge, including strictness when required.
The converse follows by testing the indicated hinges.
\end{proof}

\subsection{Uniform first-order expansion}

For a fixed rational piecewise-affine function \(f\), the expansion from the preceding section gives
\[
\mathcal C_k(f)
=
\frac{1}{kA(1)}
\left(B(f)-\frac{B(1)}{A(1)}A(f)\right)
+O_f(k^{-2})
\]
as \(k\to\infty\) through multiples of a common denominator of the endpoints and the breakpoints of \(f\). 
The remainder in the expansion depends on \(f\). 
The purpose of the next statement is to prove a control of the remainder, uniform in the hinge \(h_t\). 
This will then be applied to get a uniform positivity of the Chow weights for hinges in a compact subset of the interior, under the assumption of K-polystability. 

\begin{lem}\label{lem:rank-one-uniform-hinge-expansion}
As \(k\to\infty\), uniformly for \(t\in\Gamma_k\),
\[
\mathcal C_k(h_t)
=
\frac{1}{kA(1)}
\left(
\frac{(b-t)P(b)}{2}
+\int_t^b (x-t)Q(x)\,dx
-\frac{B(1)}{A(1)}\int_t^b (x-t)P(x)\,dx
\right)
+O(k^{-2}).
\]
\end{lem}

\begin{proof}
The hinge \(h_t\) vanishes on \([a,t]\). 
Since \(t,b\in\Gamma_k\), the trapezoidal rule on \([t,b]\) gives
\[
\sum_{j/k\in\Gamma_k} h_t(j/k)P(j/k)
=
k\int_t^b (x-t)P(x)\,dx
+\frac{(b-t)P(b)}{2}
+O(k^{-1}).
\]
The boundary contribution at \(t\) vanishes. 
The error is uniform because the interval has length at most \(b-a\) and the second derivative of \((x-t)P(x)\) is bounded independently of \(t\). 
Applying the same rule to \((x-t)Q(x)\) gives
\[
\sum_{j/k\in\Gamma_k} h_t(j/k)Q(j/k)
=
k\int_t^b (x-t)Q(x)\,dx+O(1),
\]
uniformly in \(t\), with the bounded boundary contribution included in \(O(1)\).

The leading Weyl term \(k^{n-1}P\) thus contributes \(k^nA(h_t)+k^{n-1}(b-t)P(b)/2\) up to \(O(k^{n-2})\).
The next term \(k^{n-2}Q\) contributes \(k^{n-1}\int_t^b(x-t)Q(x)\,dx\) up to the same order.
The remaining Weyl terms contribute \(O(k^{n-2})\) in total, since their coefficients and \(h_t\) are uniformly bounded and there are \(O(k)\) grid points. 
Therefore
\[
\sum_{j/k\in\Gamma_k}\dim W_{k,k\chi+j\xi}\,h_t(j/k)
=
k^nA(h_t)+k^{n-1}B(h_t)+O(k^{n-2}).
\]

Taking \(f=1\) in the proof of Proposition~\ref{prop:Chow-weight-expansion} gives the  expansion \(p_k=k^nA(1)+k^{n-1}B(1)+O(k^{n-2})\). 
Dividing the weighted sum by \(p_k\) and subtracting \(A(h_t)/A(1)\) gives
\[
\mathcal C_k(h_t)
=
\frac{1}{kA(1)}
\left(B(h_t)-\frac{B(1)}{A(1)}A(h_t)\right)
+O(k^{-2}).
\]
The error is uniform because \(A(h_t)\) and \(B(h_t)\) are bounded independently of \(t\). 
Substituting their expressions gives the stated formula.
\end{proof}

\begin{lem}\label{lem:rank-one-Chow-away-from-endpoints}
Assume that \((X,L)\) is \(G\)-equivariantly K-polystable.
Fix \(0<\varepsilon<(b-a)/2\), and set
\[
I_\varepsilon=
\begin{cases}
[a,b-\varepsilon] & \text{if \(X\) is non-horospherical},\\
[a+\varepsilon,b-\varepsilon] & \text{if \(X\) is horospherical}.
\end{cases}
\]
Then for all sufficiently large \(k\), \(\mathcal C_k(h_t)>0\) for every \(t\in\Gamma_k\cap I_\varepsilon\).
\end{lem}

\begin{proof}
By the comparison with affine test functions established in the proof of~\cite[Theorem~3.1]{Delcroix_2023_RK1},
we have
\(\mathcal L(h_t)>0\) for every real \(t\in I_\varepsilon\).
By continuity and compactness,
\[
c_\varepsilon:=\min_{t\in I_\varepsilon}\mathcal L(h_t)>0.
\]
Lemma~\ref{lem:rank-one-uniform-hinge-expansion} now gives
\[
\mathcal C_k(h_t)
\geq \frac{c_\varepsilon}{2kA(1)}-O(k^{-2})>0
\]
uniformly for \(t\in\Gamma_k\cap I_\varepsilon\) and all sufficiently large \(k\).
\end{proof}

\subsection{Control near the endpoints} 

Set \(\ell_\alpha(x)=\langle\chi+x\xi,\alpha^\vee\rangle\), and \(c_\alpha=\langle\rho,\alpha^\vee\rangle\). 
The Weyl dimension formula gives, for \(j/k\in\Gamma_k\), \(\dim W_{k,k\chi+j\xi}=W_k(j/k)\), with
\[
W_k(x)=\prod_\alpha\frac{k\ell_\alpha(x)+c_\alpha}{c_\alpha},
\]
where the product runs over the positive roots \(\alpha\) for which \(\ell_\alpha\) is not identically zero on \([a,b]\).
Let \(S_b\) be the set of roots in the Weyl product for which \(\ell_\alpha(b)=0\), and put \(m_b=|S_b|\). 
For each \(\alpha\in S_b\), set
\[
\beta_{b,\alpha}=\frac{c_\alpha}{-\ell_\alpha'},
\]
where \(\ell_\alpha'=\langle\xi,\alpha^\vee\rangle<0\) is the slope of \(\ell_\alpha\).
Let \(c_b>0\) be the leading coefficient of \(P\) at \(b\), so that, as \(y\to 0^+\), 
\[
P(b-y)=c_b y^{m_b}+O(y^{m_b+1}).
\]

We now consider hinges with breakpoint \(t=b-l/k\), where \(l\) is a positive integer and \(t\in\Gamma_k\). 
Only the grid points \(b-j/k\) with \(0\leq j<l\) contribute to the discrete term in \(\mathcal C_k(h_{b-l/k})\). 
Note that, given a fixed \(l\geq 1\), for all sufficiently large \(k\), we have \(0<l/k<b-a\), and \(b-l/k\) belongs to \(\Gamma_k\) because \(b\in\Gamma_k\).

\begin{lem}\label{lem:rank-one-endpoint-laurent}
For each integer \(l\geq 1\), there exists a Laurent polynomial \(U_{b,l}\) such that
\[
p_k\mathcal C_k(h_{b-l/k})=U_{b,l}(k)
\]
for every integer \(k\geq 1\) with \(b-l/k\in\Gamma_k\).
\end{lem}

\begin{proof}
For the discrete part, the hinge is nonzero at exactly the grid points \(b-j/k\) with \(0\leq j<l\), where its value is \((l-j)/k\). 
For fixed \(j\), the Weyl product makes \(W_k(b-j/k)\) a polynomial in \(k\), since \(
k\ell_\alpha(b-j/k)=k\ell_\alpha(b)-j\ell_\alpha'\).
By higher-cohomology vanishing for ample line bundles on
projective spherical varieties
\cite[Proposition~1.3.1(vii)]{Alexeev_Brion_2006},
\(p_k\) agrees with its Hilbert polynomial for every \(k\geq1\). 
Changing variables by \(y=k(b-x)\) in the continuous term gives
\begin{equation}\label{eq:rank-one-endpoint-laurent}
p_k\mathcal C_k(h_{b-l/k})
=
\frac{1}{k}\sum_{j=0}^{l-1}(l-j)W_k(b-j/k)
-
\frac{p_k}{A(1)k^2}\int_0^l(l-y)P(b-y/k)\,dy.
\end{equation}
The integral is a polynomial in \(k^{-1}\), since \(P\) is a polynomial and the limits of integration are fixed.
Thus the right-hand side defines the required Laurent polynomial \(U_{b,l}\).
\end{proof}

\begin{lem}\label{lem:rank-one-endpoint-expansion}
For each integer \(l\geq 1\), \(k^{m_b+2-n}U_{b,l}(k)\) is a polynomial in \(k^{-1}\) with constant term \(c_bD_b(l)\).
Here \(D_b\) is the polynomial in \(l\) whose values at positive integers are given by
\[
D_b(l)
=
\sum_{j=0}^{l-1}(l-j)
\prod_{\alpha\in S_b}(j+\beta_{b,\alpha})
-\frac{l^{m_b+2}}{(m_b+1)(m_b+2)}.
\]
\end{lem}

\begin{proof}
Multiplying~\eqref{eq:rank-one-endpoint-laurent} by \(k^{m_b+2-n}\) gives
\[
\begin{aligned}
k^{m_b+2-n}U_{b,l}(k)
={}&\sum_{j=0}^{l-1}(l-j)\,k^{m_b+1-n}W_k(b-j/k)\\
&-\frac{p_k}{A(1)k^n}
\int_0^l(l-y)\,k^{m_b}P(b-y/k)\,dy.
\end{aligned}
\]
For \(\alpha\in S_b\), the corresponding Weyl factor at \(b-j/k\) is
\[
\frac{-\ell_\alpha'}{c_\alpha}(j+\beta_{b,\alpha}).
\]
The remaining \(n-1-m_b\) factors are linear in \(k\), with leading coefficients \(\ell_\alpha(b)/c_\alpha\). 
Thus \(k^{m_b+1-n}W_k(b-j/k)\) is a polynomial in \(k^{-1}\) with constant term \(c_b\prod_{\alpha\in S_b}(j+\beta_{b,\alpha})\). 
Here \(c_b\) is the product of the factors \(-\ell_\alpha'/c_\alpha\) for \(\alpha\in S_b\) and \(\ell_\alpha(b)/c_\alpha\) for the remaining roots, as follows
from the product formula for \(P\).

In the second term, \(p_k/(A(1)k^n)\) and \(k^{m_b}P(b-y/k)\) are polynomials in \(k^{-1}\), with constant terms \(1\) and \(c_by^{m_b}\), respectively.
The integral therefore contributes the constant term
\[
c_b\int_0^l(l-y)y^{m_b}\,dy
=
\frac{c_b l^{m_b+2}}{(m_b+1)(m_b+2)}.
\]
Subtracting gives the claimed constant term \(c_bD_b(l)\).
\end{proof}

\begin{lem}\label{lem:rank-one-endpoint-polynomial-degree}
If \(m_b=0\), then \(D_b(l)=l/2\). 
If \(m_b\geq 1\), the polynomial \(D_b\) has degree \(m_b+1\) and its leading coefficient is
\[
\frac{\sum_{\alpha\in S_b}\beta_{b,\alpha}}{m_b(m_b+1)}>0.
\]
\end{lem}

\begin{proof}
If \(m_b=0\), the empty product gives \(D_b(l)=\frac{l(l+1)}{2}-\frac{l^2}{2}=\frac l2\).

Assume now that \(m=m_b\geq1\).
For \(r\geq1\), the power-sum formulas (Faulhaber's formula) give
\[
\sum_{j=0}^{l-1}j^r
=\frac{l^{r+1}}{r+1}-\frac{l^r}{2}+O(l^{r-1}).
\]
Consequently,
\[
\sum_{j=0}^{l-1}(l-j)j^r
=l\sum_{j=0}^{l-1}j^r-\sum_{j=0}^{l-1}j^{r+1} 
=\frac{l^{r+2}}{(r+1)(r+2)}+O(l^r).
\]
For \(r=0\), the corresponding sum is \(l(l+1)/2\).

Writing
\[
\prod_{\alpha\in S_b}(j+\beta_{b,\alpha})
=
j^m+
\left(\sum_{\alpha\in S_b}\beta_{b,\alpha}\right)j^{m-1}
+O(j^{m-2}),
\]
we deduce from the above power-sum formulas that
\[
\sum_{j=0}^{l-1}(l-j)
\prod_{\alpha\in S_b}(j+\beta_{b,\alpha})
=
\frac{l^{m+2}}{(m+1)(m+2)}
+
\frac{\sum_{\alpha\in S_b}\beta_{b,\alpha}}{m(m+1)}
l^{m+1}
+O(l^m).
\]
By definition of \(D_b\), the first term cancels, leaving
\[
D_b(l)=
\frac{\sum_{\alpha\in S_b}\beta_{b,\alpha}}{m(m+1)}
l^{m+1}+O(l^m).
\]
Since all \(\beta_{b,\alpha}\) are positive, this gives the claimed degree and leading coefficient.
\end{proof}

\begin{lem}\label{lem:rank-one-uniform-endpoint}
There exists a constant \(C>0\) such that, for all positive integers \(k,l\) with \(b-l/k\in\Gamma_k\),
\[
\left|k^{m_b+2-n}U_{b,l}(k)-c_bD_b(l)\right|
\leq C\frac{l^{m_b+2}}{k}.
\]
\end{lem}

\begin{proof}
Put \(m=m_b\) and \(q=n-1-m\), and write
\[
P(b-y)=\sum_{r=0}^{q}a_r y^{m+r}.
\]
In the normalized Weyl product \(k^{-q}W_k(b-j/k)\), the factors indexed by \(S_b\) are
\(1-j\ell_\alpha'/c_\alpha\). 
Each of the other factors is
\[
\frac{\ell_\alpha(b)}{c_\alpha}
+\frac{1}{k}\left(1-\frac{j\ell_\alpha'}{c_\alpha}\right).
\]
Write \(k^{-q}W_k(b-j/k)=\sum_{r=0}^{q}k^{-r}R_r(j)\).
To obtain \(R_r(j)\), one selects the \(k^{-1}\) term from exactly \(r\) of the factors outside \(S_b\). 
Each summand then has \(m+r\) factors affine in \(j\), so \(\deg R_r\leq m+r\). 
The coefficient of \(j^{m+r}\) comes from taking the linear term in each of these factors.
These same choices give the coefficient \(a_r\) of \(y^{m+r}\) in \(P(b-y)\).

Set
\[
T_r(l)=\sum_{j=0}^{l-1}(l-j)R_r(j),
\qquad
I_r(l)=\int_0^l(l-y)a_r y^{m+r}\,dy
=\frac{a_r l^{m+r+2}}{(m+r+1)(m+r+2)}.
\]
By the power-sum formulas used in the proof of Lemma~\ref{lem:rank-one-endpoint-polynomial-degree}, \(\deg T_r\leq m+r+2\), with coefficient \(a_r/((m+r+1)(m+r+2))\) at degree \(m+r+2\). 
Hence \(T_r(l)-I_r(l)=O(l^{m+r+1})\) for \(l\geq1\). 
For \(r=0\), the definition of \(D_b\) gives \(T_0(l)-I_0(l)=c_bD_b(l)\).

Substituting these expansions into formula~\eqref{eq:rank-one-endpoint-laurent} gives
\[
\begin{aligned}
k^{m+2-n}U_{b,l}(k)-c_bD_b(l)
={}&\sum_{r=1}^{q}k^{-r}\bigl(T_r(l)-I_r(l)\bigr)\\
&+\left(1-\frac{p_k}{A(1)k^n}\right)
  \sum_{r=0}^{q}k^{-r}I_r(l).
\end{aligned}
\]
Since \(l/k\leq b-a\), each term in the first sum is bounded by a constant times
\[
k^{-r}l^{m+r+1}
=\frac{l^{m+2}}{k}\left(\frac{l}{k}\right)^{r-1}.
\]
The sum involving \(I_r(l)\) is \(O(l^{m+2})\), since \(k^{-r}I_r(l)=O(l^{m+2}(l/k)^r)\). Its prefactor is \(O(k^{-1})\) by the Hilbert polynomial expansion. 
Both contributions are therefore \(O(l^{m+2}/k)\), with a constant independent of \(k\) and \(l\).
\end{proof}

\begin{cor}
There exist \(0<\varepsilon<b-a\) and an integer \(k_0\geq1\) such that, for every integer \(k\geq k_0\) and every integer \(1\leq l\leq\varepsilon k\) with \(D_b(l)>0\), we have
\[
\mathcal C_k(h_{b-l/k})>0.
\]
\end{cor}

\begin{proof}
Since \(D_b\) has degree \(m_b+1\) and positive leading coefficient, there exists \(d>0\) such that \(D_b(l)\geq d\,l^{m_b+1}\) for every positive integer \(l\) with \(D_b(l)>0\). 
Indeed, the quotient \(D_b(l)/l^{m_b+1}\) tends to a positive limit so is strictly positive from a threshold \(l_0\), and there are only finitely many integers \(l\leq l_0\). 

By Lemma~\ref{lem:rank-one-uniform-endpoint}, for \(1\leq l\leq\varepsilon k\) with \(D_b(l)>0\),
\[
k^{m_b+2-n}U_{b,l}(k)
\geq (c_b d-C\varepsilon)l^{m_b+1}.
\]
Choose \(\varepsilon\) so that \(C\varepsilon<c_b d\).
The result follows from \(U_{b,l}(k)=p_k\mathcal C_k(h_{b-l/k})\) and \(p_k>0\).
\end{proof}

If \(D_b(l)<0\) for some positive integer \(l\), then Lemma~\ref{lem:rank-one-endpoint-expansion} gives \(\mathcal C_k(h_{b-l/k})<0\) for all sufficiently large \(k\).

It remains to consider the positive integers \(l\) with \(D_b(l)=0\). 
There are only finitely many, since \(D_b\) is a nonzero polynomial. 
For each such \(l\), the Chow weight \(\mathcal C_k(h_{b-l/k})\) vanishes identically if
\(U_{b,l}\equiv0\). 
Otherwise, its sign for all sufficiently large \(k\) is the sign of the coefficient of the highest power of \(k\) in \(U_{b,l}\). 
Since the set of these \(l\) is finite, a single threshold in \(k\) suffices for all of them.

In the horospherical case, the calculation at \(a\) is similar to that at \(b\), using the reflected hinge \((t-x)_+\) with \(t=a+l/k\). 
Define \(D_a\) by the same formula as \(D_b\), using \(S_a=\{\alpha:\ell_\alpha(a)=0\}\), \(m_a=|S_a|\), and \(\beta_{a,\alpha}=c_\alpha/\ell_\alpha'\).
Let \(U_{a,l}\) be the Laurent polynomial giving \(p_k\mathcal C_k((a+l/k-x)_+)\). 
Since \(h_t(x)-(t-x)_+=x-t\), the reflected hinge has the same Chow weight as \(h_t\) whenever \(\mathcal C_k(x)=0\).

\subsection{A finite test for asymptotic Chow polystability}

\begin{thm}\label{thm:rank-one-asymptotic-Chow}
Assume one of the following:
\begin{enumerate}[label=(\roman*)]
\item \(X\) is non-horospherical and \((X,L)\) is \(G\)-equivariantly K-polystable;
\item \(X\) is horospherical and \(\mathcal C_k(x)=0\) at \(n=\dim X\) distinct positive integers \(k\).
\end{enumerate}
Set \(E=\{b\}\) in case~{(i)} and \(E=\{a,b\}\) in case ~{(ii)}.

Then \((X,L)\) is asymptotically Chow semistable if and only if the following conditions hold for every \(e\in E\):
\begin{enumerate}
\item \(D_e(l)\geq0\) for every positive integer \(l\).
\item At every positive integer zero \(l\) of \(D_e\), the Laurent polynomial \(U_{e,l}\) is identically zero or has positive leading coefficient.
\end{enumerate}
It is asymptotically Chow polystable if and only if the first condition holds and every Laurent polynomial in the second condition has positive leading coefficient.
\end{thm}

\begin{proof}
Recall that \(G\)-equivariant K-polystability is equivalent to \(\mathcal L(x)>0\) in the non-horospherical case and to \(\mathcal L(x)=0\) in the horospherical case, see~\cite[Theorem~3.1 and subsequent discussion]{Delcroix_2023_RK1}.
In particular, in the horospherical case the finite vanishing condition in Theorem~\ref{thm:rank-one-asymptotic-Chow}  implies \(G\)-equivariant K-polystability. 
Indeed, Proposition~\ref{prop:product-Chow-weight-finite-test} gives \(\mathcal C_k(x)=0\) for all \(k\geq1\), and the first-order expansion of \(\mathcal C_k(x)\) then gives \(\mathcal L(x)=0\).
In the horospherical case, this also shows that the Chow weight of the reflected hinge at \(a+l/k\) equals \(\mathcal C_k(h_{a+l/k})\).

For necessity, fix \(e\in E\) and an integer \(l\geq1\).
The breakpoint \(b-l/k\) if \(e=b\), or \(a+l/k\) if \(e=a\), belongs to \(\Gamma_k\cap(a,b)\) for all sufficiently large \(k\). 
If \(D_e(l)<0\), the endpoint expansion gives a negative Chow weight at this breakpoint for all sufficiently large \(k\). 
If \(D_e(l)=0\), the eventual sign of the Chow weight is the sign of the leading coefficient of \(U_{e,l}\) when \(U_{e,l}\) is nonzero, and the weight is zero when
\(U_{e,l}\equiv0\). 
Proposition~\ref{prop:rank-one-finite-Chow-tests} therefore gives the necessary conditions for semistability and polystability.

For sufficiency, choose \(0<\varepsilon<(b-a)/2\) small enough that the uniform positivity result near the endpoints applies at every \(e\in E\). 
For all sufficiently large \(k\), Lemma~\ref{lem:rank-one-Chow-away-from-endpoints} gives strict positivity at every grid breakpoint in \(I_\varepsilon\).
Every remaining breakpoint required by Proposition~\ref{prop:rank-one-finite-Chow-tests} is at distance \(l/k<\varepsilon\) from some \(e\in E\), with \(l\geq1\).
If \(D_e(l)>0\), the uniform endpoint estimate gives strict positivity at that breakpoint. 
If \(D_e(l)=0\), the condition on \(U_{e,l}\) gives non-negativity for all sufficiently large \(k\), and strict positivity under the polystability condition.
Since there are only finitely many pairs \((e,l)\) with \(D_e(l)=0\), one threshold in \(k\) suffices for all of them. 
The finite Chow criterion now gives asymptotic Chow semistability or polystability, respectively.
\end{proof}

\begin{rem}
The conditions in Theorem~\ref{thm:rank-one-asymptotic-Chow} can be checked at finitely many integers. 
If \(m_e=0\), then \(D_e(l)=l/2>0\) for every \(l\geq1\). 
If \(m_e\geq1\), the power-sum formulas give \(D_e(0)=0\), so \(q_e(l)=D_e(l)/l\) is a polynomial of degree \(m_e\) with positive leading coefficient.

The polynomial \(q_e\) is monotone between consecutive critical points and tends to \(+\infty\) as \(l\to+\infty\).
Its minimum over the positive integers is therefore attained at \(1\) or at one of \(\lfloor t\rfloor,\lceil t\rceil\) for a real critical point \(t\geq1\). 
It suffices to check \(D_e(l)\geq0\) at those integers. 
Since \(q_e'\) has degree \(m_e-1\), there are at most \(2m_e-1\) integer candidates, including \(1\).
Moreover, if \(m_e\geq2\) and \(q_e(t)=\sum_{j=0}^{m_e}a_jt^j\), the Cauchy root bound
applied to \(q_e'\) gives 
\[
|t|\leq 1+\max_{1\leq j<m_e}\frac{j|a_j|}{m_e a_{m_e}}
\]
for every real critical point \(t\) of \(q_e\).
The Laurent polynomial test is then needed only at the positive integer roots of \(D_e\) which similarly are in a finite set.
\end{rem}

We now give a simple sufficient condition when \(P\) vanishes to order at most one at the relevant endpoint. 

\begin{cor}
Under the hypotheses of Theorem~\ref{thm:rank-one-asymptotic-Chow}, an endpoint \(e\in E\) with \(P(e)>0\) satisfies \(D_e(l)>0\) for every positive integer \(l\). 
If \(P\) has a simple zero at \(e\), then \(D_e(l)>0\) for every positive integer \(l\) if and only if
\[
6Q(e)>|P'(e)|.
\]
If \(6Q(e)=|P'(e)|\), the only positive integer zero of \(D_e\) is \(l=1\), and \(D_e(l)>0\) for \(l\geq2\).
\end{cor}

\begin{proof}
If \(P(e)>0\), then \(m_e=0\) and \(D_e(l)=l/2\).
If \(P\) has a simple zero at \(e\), exactly one Weyl factor vanishes there. 
Writing \(\beta\) for its shift \(\beta_{e,\alpha}\), the Weyl product gives \(Q(e)=\beta|P'(e)|\). 
The defining sum becomes
\[
D_e(l)
=
\sum_{j=0}^{l-1}(l-j)(j+\beta)-\frac{l^3}{6}
=
\frac{l}{6}\bigl(3\beta(l+1)-1\bigr).
\]
Since \(\beta>0\), the factor in parentheses is strictly increasing in \(l\). 
Its value at \(l=1\) is \(6\beta-1\), which gives both assertions.
\end{proof}

\subsection{A K-polystable non-toric variety that is  
asymptotically Chow unstable
}

\label{exa:critical-endpoint-quotient}
Let \(G=\mathrm{SL}_3\), and denote by \(\omega_1,\omega_2\) its fundamental weights.
Consider the polarized horospherical \(G\)-variety \((X,L)\) with weight lattice \(M=\mathbb Z\bigl(6(\omega_1-\omega_2)\bigr)\) and moment polytope \(\Delta=\operatorname{conv}(6\omega_1,6\omega_2)\). 
The parabolic subgroup normalizing the horospherical subgroup is the Borel \(B\), so that \(X\) has dimension four and spherical rank one.
The section decomposition is
\[
H^0(X,L^{\otimes k})
=
\bigoplus_{j=0}^k
W_{k,6j\omega_1+6(k-j)\omega_2},
\]
with multiplication given by the Cartan product.
In particular, \(L\) is very ample.

Geometrically, \(X\) can be realized as a quotient of the smooth quadric
\[
Q^4=\{[v:\varphi]\in
\mathbb P(\mathbb C^3\oplus(\mathbb C^3)^*):
\varphi(v)=0\}
\]
by the action of \(\mu_6\) given by \(\zeta\cdot[v:\varphi]=[\zeta v:\varphi]\).
The variety \(X\) is not toric. 
Indeed, a normal finite cover of a toric variety which is \'etale in codimension one is itself toric.
If \(X\) were toric, then \(Q^4\) would thus be a smooth projective toric variety of Picard number one, hence isomorphic to \(\mathbb P^4\), a contradiction.

The singular locus of \(X\) consists of the images of the two disjoint projective planes \(\{v=0\}\) and \(\{\varphi=0\}\).
Along each of them, the transverse singularities are cyclic quotient singularities of type \(\frac16(1,1)\).
In particular, \(X\) has klt but noncanonical singularities.

We take \(\chi=6\omega_2\),
\(\xi=6(\omega_1-\omega_2)\), and \([a,b]=[0,1]\).
The Weyl dimension formula gives
\[
W_k(x)
=(3k+1)(6kx+1)(6k(1-x)+1),
\]
and consequently \(P(x)=108x(1-x)\), \(Q(x)=36x(1-x)+18\).
By symmetry, \(\mathcal C_k(x)=0\) for every \(k\geq1\).

The vanishing of \(\mathcal C_k(x)\) implies \(\mathcal L(x)=0\), hence \(G\)-equivariant K-polystability, as observed in the proof of Theorem~\ref{thm:rank-one-asymptotic-Chow}.
Since \(X\) is \(\mathbb Q\)-Fano and \(L\sim_{\mathbb Q}-\frac32 K_X\), the equivariant reduction theorem~\cite{Zhuang_2021} gives K-polystability of \((X,L)\).

We now consider asymptotic Chow stability.
At each endpoint, exactly one Weyl factor vanishes to leading order, with \(\beta_{e,\alpha}=1/6\).
Consequently,
\[
D_a(l)=D_b(l)=\frac{l(l-1)}{12}.
\]
These polynomials are strictly positive for every integer \(l\geq2\), but \(D_a(1)=D_b(1)=0\).

To resolve these zeros, it suffices by symmetry to compute \(U_{b,1}\). 
Since \(p_k=(3k+1)(k+1)(6k^2+1)\), formula~\eqref{eq:rank-one-endpoint-laurent} gives
\[
\begin{aligned}
U_{b,1}(k)
&=\frac{(3k+1)(6k+1)}{k}
-p_k\left(\frac1{k^3}-\frac1{2k^4}\right)\\
&=-\frac{(3k+1)(k-1)(4k^2+1)}{2k^4}.
\end{aligned}
\]
In particular, the leading coefficient of \(U_{b,1}\) is negative, and Theorem~\ref{thm:rank-one-asymptotic-Chow} gives asymptotic Chow instability.
More precisely, \(U_{b,1}(k)<0\) for every \(k\geq2\), so every embedding defined by \(L^{\otimes k}\) with \(k\geq2\) is Chow unstable.
In contrast, the embedding defined by \(L\) is Chow polystable, since \(\mathcal C_1(x)=0\) and there are no interior grid hinges.

Thus, even for non-toric horospherical varieties of rank one, K-polystability together with the vanishing of the Chow weights of all \(G\)-equivariant product test configurations at every level does not imply asymptotic Chow semistability.

\subsection{Chow stability of a wonderful \(G_2\)-variety}
\label{subsec:G2-Chow}

Consider the smooth wonderful \(G_2\)-variety \(X\) of~\cite[\S\S4.1.2--4.1.3]{Delcroix_2023_RK1}. 
Let \(P_1\subset G_2\) be the maximal parabolic subgroup associated with the short simple root.
The homogeneous space \(G_2/P_1\) is the smooth five-dimensional quadric \(Q^5\).
Consider the associated bundle 
\[
X=G_2\times^{P_1}
\mathbb P(\operatorname{Sym}^2\mathbb C^2),
\]
where \(P_1\) acts on \(\operatorname{Sym}^2\mathbb C^2\) through its Levi quotient \(\operatorname{GL}_2(\mathbb{C})\). 
The fiber is thought of as the space of conics in \(\mathbb{P}^1\). 
Thus \(X\) is a smooth sevenfold, fibered in projective planes over \(Q^5\).
It is a wonderful \(G_2\)-variety with a unique closed orbit \(E\simeq G_2/B\), corresponding fiberwise to degenerate quadratic forms.

Let \(\pi:X\to G_2/P_1\simeq Q^5\) denote the natural projection. 
The color \(D_1\) is the pullback of a Schubert divisor on \(Q^5\), whereas \(D_2\) restricts on each fiber \(\mathbb P^2\) to the hyperplane of conics containing a fixed point of \(\mathbb P^1\).
We refer to~\cite[\S4.1.2]{Delcroix_2023_RK1} for further details. 
The classes of the colors \(D_1,D_2\) form a \(\mathbb Z\)-basis of \(\operatorname{Pic}(X)\).
The ample line bundles are precisely
\[
L_{u,v}=\mathcal O_X(uD_1+vD_2),
\]
for \(u,v\in\mathbb Z_{>0}\), and they are all very ample~\cite[Theorem 15]{Pezzini_2007}.

Number the simple roots so that \(\alpha_1\) is short, and denote the fundamental weights by \(\omega_1,\omega_2\). 
The spherical root is \(\xi=2\alpha_2\).
The moment polytope for \(L_{u,v}\) is
\[
\Delta_{u,v}
=\chi+[-v/2,0]\xi,
\qquad
\chi=u\omega_1+2v\omega_2.
\]
The decomposition of sections is thus 
\[
H^0(X,L_{u,v}^{\otimes k})
=
\bigoplus_{i=0}^{\lfloor kv/2\rfloor}
W_{k,k\chi-i\xi}.
\]
In particular, \(\Gamma_k=[-v/2,0]\cap k^{-1}\mathbb Z\).

The Weyl dimension formula gives
\[
W_k(x)=\frac1{120}
\begin{aligned}[t]
&\bigl(k(u-6x)+1\bigr)
 \bigl(k(2v+4x)+1\bigr)\\
&\cdot\bigl(k(u+2v-2x)+2\bigr)
 \bigl(k(u+4v+2x)+3\bigr)\\
&\cdot\bigl(k(u+6v+6x)+4\bigr)
 \bigl(k(2u+6v)+5\bigr).
\end{aligned}
\]

Following~\cite{Delcroix_2023_RK1}, we introduce the parameter 
\[
s=\frac{3v}{u+3v}\in(0,1)
\]
which parametrizes ample line bundles up to multiple. 
We will compare K-polystability and Chow stability throughout this family of polarizations.

Recall from~\cite[\S4.1.3]{Delcroix_2023_RK1} that \(L_{u,v}\) is K-polystable precisely when \(f(s)>0\), where
\[
f(s)=11s^8+20s^7-348s^6-240s^5+3123s^4+1260s^3-9072s^2+5103.
\]
The polynomial \(f\) has a unique root \(s_0\in(0,1)\), with \(s_0\simeq0.972020421341\).
Thus \(L_{u,v}\) is K-polystable if and only if \(s<s_0\).

The root \(s_0\) is irrational. 
Indeed, suppose that \(p/q\) is a rational root of \(f\), with \(p,q\in\mathbb Z\) coprime and \(q>0\). 
By the rational root theorem, \(q\mid 11\).
In particular, \(q\) would be odd. 
Reducing \(q^8f(p/q)=0\) modulo \(2\) would then yield \(p^8+p^4\equiv 1\pmod 2\), which is impossible. 

Consequently, every ample rational polarization is either K-polystable or K-unstable.

Since \(X\) is non-horospherical, only the right endpoint \(b=0\) enters Theorem~\ref{thm:rank-one-asymptotic-Chow}.
This endpoint corresponds to the strictly dominant weight \(\chi=u\omega_1+2v\omega_2\), so \(P(b)>0\).
Consequently, \(D_b(l)=l/2>0\) for every integer \(l\geq1\).

Thus K-polystability implies asymptotic Chow polystability. 
Conversely, asymptotic Chow semistability implies K-semistability, which is equivalent to K-polystability for rational polarizations by the irrationality of \(s_0\).
We conclude that
\[
\begin{aligned}
(X,L_{u,v})\text{ is K-polystable}
&\iff (X,L_{u,v})\text{ is asymptotically Chow polystable}\\
&\iff (X,L_{u,v})\text{ is asymptotically Chow semistable}.
\end{aligned}
\]

On the other hand, fixed level embeddings are often Chow unstable. 

\begin{exa}
Consider the polarization \(L=L_{1,9}\), for which \(s=\frac{27}{28}<s_0\).
Thus \((X,L)\) is K-polystable and asymptotically Chow polystable. 
Nevertheless, its embedding defined by \(L\) is Chow unstable.
Indeed, the degree-one weight grid is
\(
\Gamma_1=\{-4,-3,-2,-1,0\}.
\)
Direct computation using the Weyl dimension formula gives
\[
\mathcal C_1(x)
=\frac{648273693}{402128218492}>0,
\]
whereas, for the hinge \(h_{-4}(x)=\max\{0,x+4\}\),
\[
\mathcal C_1(h_{-4})
=-\frac{231059245295}{97717157093556}<0.
\]
By Proposition~\ref{prop:rank-one-finite-Chow-tests},
the embedding is Chow unstable.
\end{exa}

This example shows that positivity of the affine Chow test, even for a K-polystable polarization, does not imply Chow semistability (unlike for K-polystability).

\begin{prop}\label{prop:G2-no-uniform-Chow-exponent}
For every integer \(k_0\geq1\), there are infinitely many primitive ample line bundles \(L\) on \(X\) such that \((X,L)\) is K-polystable and asymptotically Chow polystable, but the embedding defined by \(L^{\otimes k_0}\) is Chow unstable.

In particular, there is no uniform exponent \(k_*\) such that every K-polystable ample line bundle on \(X\) is Chow polystable at all levels \(k\geq k_*\).
\end{prop}

\begin{proof}
Fix \(k_0\geq1\), and set \(C=\frac{u+3v}{6}\), \(K=k_0C\).
In the coordinate \(t=(x+v/2)/C\), the moment segment becomes \([0,s]\), and the occurring weights at level \(k_0\) correspond to \(t=\frac jK\), for \(j=Ks-i\) and \(0\leq i\leq\lfloor Ks\rfloor\).
Since \(Ks=k_0v/2\), the values of \(j\) are integers if \(k_0v\) is even and half-integers otherwise.
The Weyl dimension formula gives
\[
\begin{aligned}
w_K(j)=\frac1{60}
&(6K-6j+1)(4j+1)(6K+6j+4)\\
&\cdot(12K+5)(3K-j+1)(6K+2j+3).
\end{aligned}
\]
Write
\[
P_*(t)=\frac{288}{5}t(1-t^2)(9-t^2),
\qquad
A_i(s)=\int_0^s t^iP_*(t)\,dt.
\]
Since \(x=Ct-v/2\) and \(\mathcal C_{k_0}\) vanishes on constants, we have
\[
\mathcal C_{k_0}(x)
=
\frac{C}{A_0(s)\sum_j w_K(j)}N(K,s),
\]
where
\[
N(K,s)=
A_0(s)\sum_j\frac jK w_K(j)
-A_1(s)\sum_j w_K(j).
\]
Thus \(\mathcal C_{k_0}(x)\) and \(N(K,s)\) have the same sign.

Expanding the Weyl dimension polynomial, we obtain
\[
w_K(Kt)=K^6P_*(t)+K^5Q_*(t)+K^4R_*(t)+O(K^3),
\]
where
\[
\begin{aligned}
Q_*(t)&=\frac{24}{5}(5t^5+15t^4-150t^3-90t^2+225t+27),\\
R_*(t)&=6(t^2-10t-3)(5t^2+6t-15).
\end{aligned}
\]
When \(k_0v\) is even, the Euler--Maclaurin formula gives, for \(i=0,1\),
\[
\sum_j\left(\frac jK\right)^i w_K(j)
=K^7A_i(s)+K^6B_i(s)+K^5C_i(s)+O(K^4),
\]
where, using \(P_*(0)=0\),
\[
B_i(s)=\int_0^s t^iQ_*(t)\,dt+\frac{s^iP_*(s)}2
\]
and
\[
C_i(s)=\int_0^s t^iR_*(t)\,dt 
+\frac12\left(t^iQ_*(t)\big|_{t=0}
+s^iQ_*(s)\right)
+\frac1{12}\left[(t^iP_*(t))'\right]_0^s.
\]

When \(k_0v\) is odd, the grid is shifted by half a step at the lower endpoint. 
Since \(P_*(0)=0\), the coefficient \(B_i(s)\) is unchanged. 
The lower endpoint contribution to \(C_i(s)\) becomes
\(
\frac1{24}(t^iP_*(t))'\big|_{t=0}
\)
instead of
\(
\frac12t^iQ_*(t)\big|_{t=0}
-\frac1{12}(t^iP_*(t))'\big|_{t=0}.
\)
These expressions coincide for \(i=1\) since \(P_*(0)=0\), and for \(i=0\) since \(Q_*(0)=P_*'(0)/4\).
Thus the expansion holds for both parities, with a remainder uniform for \(s\) near \(s_0\).

Substituting into the definition of \(N(K,s)\), we obtain
\[
N(K,s)=K^6(A_0B_1-A_1B_0)
+K^5(A_0C_1-A_1C_0)+O(K^4).
\]
Direct computation of these polynomials gives
\[
N(K,s)=
\frac{576}{175}K^6s^4f(s)
+\frac{48}{175}K^5s^3H(s)
+O(K^4),
\]
where
\[
\begin{aligned}
H(s)={}&55s^9+600s^8-2490s^7-14940s^6\\
&+37125s^5+107244s^4-141750s^3\\
&-183708s^2+127575s+61236.
\end{aligned}
\]

Finally, \(H(s_0)<0\). 
Indeed,
\(
0.9720204213409<s_0<0.9720204213410,
\)
and substitution of these rational bounds gives \(H(s_0)<-4683\).

Put
\(
r_0=\frac{3(1-s_0)}{s_0}.
\)
Since \(r_0\) is irrational, continued-fraction convergents from above provide infinitely many coprime positive integers \(u,v\), with \(v\to\infty\), such that
\(
0<\frac uv-r_0=O(v^{-2}).
\)
For these pairs, we have \(s<s_0\) and
\[
s_0-s
=
\frac{3(u/v-r_0)}{(r_0+3)(u/v+3)}
=O(v^{-2}).
\]
Since \(f(s_0)=0\), it follows that \(f(s)=O(v^{-2})\).
Using \(Ks=k_0v/2\), the expansion above gives
\[
\frac{175N(K,s)}{48K^5s^3}
=6k_0vf(s)+H(s)+O(K^{-1})
\longrightarrow H(s_0)<0.
\]
Thus \(\mathcal C_{k_0}(x)<0\) for all sufficiently large
\(v\), proving Chow instability.

On the other hand, \(s<s_0\) implies K-polystability and asymptotic Chow polystability, as established above. 
Finally, \(\gcd(u,v)=1\) ensures that \(L_{u,v}\) is primitive in \(\operatorname{Pic}(X)\).
\end{proof}

\section{Degenerations and faces of the Chow polytope}
\label{sec:degen}

We return to the geometric interpretation of the equivariant Chow polytope.
We describe the limit cycles associated with one-parameter subgroups of the centralizer of \(G\), relating their irreducible components to the satellites of the open orbit and to the subdivision of the moment polytope induced by the spherical envelope. 
This yields a description of the faces of the equivariant Chow polytope as Minkowski sums of smaller Chow polytopes. 
In the horospherical case, we recover an analogue of the classical description of secondary polytopes by GKZ vectors.

\subsection{Satellites and partial horospherical contractions}

We introduce the notions of satellites of a spherical subgroup \(H\) from~\cite{Batyrev_Moreau_2020} and of partial horospherical contractions from~\cite{Kaveh_Manon_2019}. 

For \(v\in\mathcal V\), let
\[
I(v)
=
\{\alpha\in\Sigma\mid\langle\alpha,v\rangle=0\}.
\]
The subset \(I(v)\) determines the smallest face of \(\mathcal V\) containing \(v\) in its interior.

For a nonzero primitive lattice point \(v\in\mathcal V\cap N\), let \(X_v\) be the elementary toroidal embedding of \(G/H\) corresponding to the ray \(\mathbb Q_{\geq0}v\), and let \(D_v\) be its closed \(G\)-orbit. 
The complement of the zero section in the normal bundle of \(D_v\) in \(X_v\) is a single \(G\)-orbit.
Let \(H_v\) be the stabilizer of a point in this orbit, and we treat \(v=0\) separately by setting \(H_0=H\). 
By~\cite[Theorem~1.1]{Batyrev_Moreau_2020}, the conjugacy class of \(H_v\) depends only on \(I(v)\). 
For \(I\subseteq\Sigma\), we denote by \(H_I\) the resulting subgroup for any such \(v\) with \(I(v)=I\), and call it the \emph{satellite} of \(H\) associated with \(I\).
The homogeneous spherical space \(G/H_I\) has the same weight lattice \(M\) as \(G/H\), and its set of spherical roots is \(I\). 
In particular, its valuation cone is
\[
\mathcal V_I
=
\left\{
w\in N_\mathbb{R}
\;\middle|\;
\langle\alpha,w\rangle\leq0
\text{ for all }\alpha\in I
\right\}.
\]
If \(I=I(v)\), then \(\mathcal V_I=\mathcal V+\mathbb{R} v\).
These properties follow from
\cite[Theorem~4.4 and Corollary~4.7]{Batyrev_Moreau_2020}.

We have \(I(v)=\Sigma\) if and only if \(v\) is in the linear part \(\mathcal{V}\cap -\mathcal{V}\) of the valuation cone, and this is further equivalent to the satellite \(H_{I(v)}\) being conjugate to \(H\). 
At the other extreme, \(H_\emptyset\) is
horospherical, and \(I(v)=\emptyset\) if and only if \(v\) is in the relative interior of \(\mathcal{V}\). 

Let
\[
R=\bigoplus_{k\geq0}H^0(X,L^{\otimes k})
 =\bigoplus_{k\geq0}
   \bigoplus_{\lambda\in\Lambda_k}W_{k,\lambda}
\]
be the section ring.
Fix a weight \(\chi\) such that
\(\Lambda_k\subset k\chi+M\).
For an integral \(v\in\mathcal V\cap N\), consider the decreasing filtration
\[
F_v^p R_k
=
\bigoplus_{\substack{\lambda\in\Lambda_k\\
\langle\lambda-k\chi,v\rangle\geq p}}
W_{k,\lambda}.
\]
It is compatible with multiplication, by the description of the valuation cone and of the multiplication of sections~\eqref{eq:section-ring-multiplication}.
The associated graded algebra
\[
R^{(v)}=\operatorname{gr}_{F_v}R
\]
has the same decomposition into irreducible \(G\)-modules as \(R\), but its multiplication is modified.
More precisely, a nonzero component
\( W_{k,\lambda}\cdot W_{\ell,\mu}
\longrightarrow
W_{k+\ell,\lambda+\mu-\gamma} \), for some \( \gamma\in\mathcal V^\vee\) survives in \(R^{(v)}\) if and only if \( \langle\gamma,v\rangle=0\). 
Thus the surviving tails belong to
\(
\mathcal V^\vee\cap v^\perp
=
\operatorname{cone}(I(v)).
\)
In particular, the resulting algebra depends only on \(I(v)\).
When \(v\) belongs to the relative interior of \(\mathcal V\), all nonzero tails disappear, and one recovers the usual horospherical contraction~\cite{Popov_1986}.
This construction is called a \emph{partial horospherical contraction} in~\cite{Kaveh_Manon_2019}.

\begin{prop}\label{prop:partial-contraction-satellite}
Let \((X,L)\) be a polarized spherical \(G\)-variety with open orbit \(G/H\), and let \(v\in\mathcal V\cap N\). 
The partial horospherical contraction
\(
X^{(v)}:=\operatorname{Proj}R^{(v)}
\)
is a projective spherical \(G\)-variety. 
Its open \(G\)-orbit is isomorphic to \(G/H_{I(v)}\).
\end{prop}

\begin{proof}
For \(v=0\), the statement is immediate.
Replacing a nonzero \(v\) by the primitive generator of its ray only reindexes the filtration, so we may assume that \(v\) is primitive.

Consider the Rees algebra
\(
\mathcal R_v
=
\bigoplus_{p\in\mathbb Z}\tau^{-p}F_v^pR
\subset R[\tau,\tau^{-1}].
\)
By the description of spherical test configurations in~\cite[Theorem~4.1 and Sections~4.2--4.3]{Delcroix_2023_KSSV2}, the family
\(
\mathcal X^{(v)}
=
\operatorname{Proj}_{\mathbb C[\tau]}\mathcal R_v
\)
is the special test configuration associated with the integral affine function \(x\mapsto\langle x-\chi,v\rangle\), to which we may add a sufficiently large integer constant without changing the underlying family.
In particular, the family
\(
\mathcal X^{(v)}
=
\operatorname{Proj}_{\mathbb C[\tau]}\mathcal R_v
\)
is normal and its central fiber
\(
\operatorname{Proj}(\mathcal R_v/(\tau))
=
\operatorname{Proj}R^{(v)}
=
X^{(v)}
\)
is normal.

As in the proof of Lemma~\ref{lem:finite-map-normalization-rees},
we have
\(
\mathcal R_v
=
R[\tau,\tau^{-1}]\cap\mathcal O_{\widehat\nu},
\)
where \(\widehat\nu\) is the invariant valuation determined by
\[
\widehat\nu(\tau^a s)
=a+\langle\lambda-k\chi,v\rangle,
\qquad 0\neq s\in W_{k,\lambda}.
\]
Since \(\widehat\nu(\tau)=1\), the ideal \((\tau)\) in \(\mathcal R_v\) is the prime ideal of elements of strictly positive valuation.
The restriction of \(\widehat\nu\) to the weight lattice of \(G/H\times\mathbb C^*\) is represented by \((v,1)\), where the last coordinate records the value of the valuation on \(\tau\). 
Thus \(\mathcal X^{(v)}\) contains the elementary toroidal embedding corresponding to the ray generated by \((v,1)\), whose closed orbit is the open \(G\times\mathbb C^*\)-orbit of \(X^{(v)}\).

On the other hand, by the proof of \cite[Proposition~2.4]{Batyrev_Moreau_2020}, the deformation of \(X_v\) to the normal bundle of \(D_v\) contains the elementary toroidal embedding corresponding to the ray generated by \((v,1)\).
Its closed orbit is the complement of the zero section in the normal bundle, hence is \(G\)-equivariantly isomorphic to \(G/H_v\).

By the Luna--Vust correspondence, these two elementary embeddings are isomorphic.
Consequently, \(X^{(v)}\) has open \(G\)-orbit \(G/H_v\simeq G/H_{I(v)}\), which concludes the proof.
\end{proof}

\subsection{Decomposition of the limit cycle}

We retain the setup of Theorem~\ref{thm:spherical-secondary-Chow}.
Let \(\eta:\mathbb C^*\to\widetilde T_0\) be a one-parameter subgroup, and let \(g=g_h\) be the spherical envelope associated with the heights \(h_\lambda=\langle\varepsilon_\lambda,\eta\rangle\).
After the finite base change used in Lemma~\ref{lem:finite-map-normalization-rees}, the source of the finite morphism \(\Pi:\mathcal X_h\to\mathcal Y_h\) has special fiber
\(
\mathcal X_0=\operatorname{Proj}(\operatorname{gr}_g R).
\)
We describe the components of \(\mathcal X_0\) and their contributions to the limiting Chow cycle
\(
\lim_{\tau\to0}\eta(\tau)^{-1}\cdot\pi_*[X].
\)

For each maximal cell \(\sigma\) of the subdivision of \(\Delta\) induced by \(g\), let \(v_\sigma\in\mathcal V\) be the slope of \(g\) on \(\sigma\) and set
\[
I_\sigma=\{\alpha\in\Sigma:\langle\alpha,v_\sigma\rangle=0\}.
\]

\begin{thm}\label{thm:limit-cycle-decomposition}
The special fiber \(\mathcal X_0\) is reduced, and its irreducible components \(X_\sigma\) are indexed by the maximal cells \(\sigma\).
Each \(X_\sigma\) is normal, with open orbit \(G/H_{I_\sigma}\)
and moment polytope \(\sigma\) for the polarization
\(\pi_\sigma^*\mathcal O_{\mathbb P(E^*)}(1)\). 
Let \(Y_\sigma=\Pi_0(X_\sigma)\) and write \(\pi_\sigma:X_\sigma\to Y_\sigma\) for the induced finite morphism. 
Then
\[
\lim_{\tau\to0}\eta(\tau)^{-1}\cdot\pi_*[X]
=\sum_\sigma(\pi_\sigma)_*[X_\sigma]
=\sum_\sigma
[\mathbb C(X_\sigma):\mathbb C(Y_\sigma)]\,[Y_\sigma].
\]
\end{thm}

\begin{proof}
By the proof of Lemma~\ref{lem:finite-map-normalization-rees}, \(\mathcal R_h\) is the intersection of \(R[\tau,\tau^{-1}]\) with the valuation rings corresponding to the affine pieces of \(g\).
Each valuation takes the value \(1\) on \(\tau\), so \((\tau)\) is the intersection of their positive-value ideals in \(\mathcal R_h\).

Let \(\ell_\sigma\) be the affine extension of \(g|_\sigma\).
The corresponding valuation takes the value \(k(\ell_\sigma-g)(\mu/k)\) on each nonzero element of \(\tau^{-kg(\mu/k)}W_{k,\mu}\). 
By the minimum formula for these valuations, the corresponding prime in \(\operatorname{gr}_g R=\mathcal R_h/(\tau)\) is the sum of the summands with \(\mu/k\notin\sigma\). 
Weights with \(\mu/k\) in the relative interiors of the maximal cells show that these primes are pairwise incomparable. 
They are therefore precisely the minimal primes, and \(\mathcal X_0\) is reduced.

Write \(A_\sigma\) for the quotient by the prime corresponding to \(\sigma\), so that \(X_\sigma=\operatorname{Proj}A_\sigma\).
Its weight monoid for \(G\times\mathbb C^*\) is the weight monoid of \(R\) intersected with the cone over \(\sigma\), and is therefore saturated.
By Brion's normality criterion~\cite[Theorem~2.14]{Brion_CIRM_2010}, \(A_\sigma\) is normal, and hence so is \(X_\sigma\).
For \(k\gg0\), its degree-\(k\) part gives the sections of \(\pi_\sigma^*\mathcal O_{\mathbb P(E^*)}(k)\).
Thus the moment polytope of \(X_\sigma\) is \(\sigma\).

Let \(\nu_\sigma\) be the valuation corresponding to \(\ell_\sigma\). 
Its restriction to the function field of \(\mathcal X_h\) is the divisorial valuation associated with \(X_\sigma\), since its center is \(X_\sigma\) and \(\nu_\sigma(\tau)=1\).
Taking quotients of highest-weight sections of the same degree, as in the proof of Proposition~\ref{prop:partial-contraction-satellite}, we see that this valuation is represented by \((v_\sigma,1)\).
Thus the corresponding elementary toroidal embedding of \(G/H\times\mathbb C^*\) is the same as in the partial horospherical contraction associated with \(v_\sigma\).
In particular, the open orbit of \(X_\sigma\) is \(G/H_{I_\sigma}\).

Finally, by Lemma~\ref{lem:finite-pushforward-Chow-limit}, the limiting cycle is the associated \(n\)-cycle of \((\Pi_0)_*\mathcal O_{\mathcal X_0}\), which is \((\Pi_0)_*[\mathcal X_0]\). 
Since \(\mathcal X_0\) is reduced, this is \(\sum_\sigma(\pi_\sigma)_*[X_\sigma]\). 
Each finite pushforward is
\[
(\pi_\sigma)_*[X_\sigma]
=
[\mathbb C(X_\sigma):\mathbb C(Y_\sigma)]\,[Y_\sigma],
\]
which gives the claimed formula.
\end{proof}

By~\cite[Lemma~3.1.3 and Lemma~2.2.2]{Alexeev_Brion_2006}, the special fiber \(\mathcal X_0\), equipped with the pullback polarization, is a polarized stable spherical variety.
The above theorem makes explicit its component moment polytopes, in the sense of~\cite[Lemma~2.3.1]{Alexeev_Brion_2006}, and identifies the open-orbit stabilizers of its irreducible components in terms of satellites of \(H\).

\subsection{Exposed faces of the Chow polytope}

Put \(Z=\pi_*[X]\) and \(Z_\sigma=(\pi_\sigma)_*[X_\sigma]\).
For a polytope \(P\), write \(\operatorname{Face}_\eta(P)\) for the face of \(P\) on which pairing with \(\eta\) is maximal:
\[
\operatorname{Face}_\eta(P)
=
\{\xi\in P:
\langle\xi,\eta\rangle
=\max_P\langle\cdot,\eta\rangle\}.
\]

\begin{cor}\label{cor:Chow-face-decomposition}
One has
\[
\operatorname{Face}_\eta
\bigl(\operatorname{Chow}_{\mathrm{GL}}^G(Z)\bigr)
=
\sum_\sigma \operatorname{Chow}_{\mathrm{GL}}^G(Z_\sigma),
\]
where the sum is a Minkowski sum in \(X^*(\widetilde T_0)_{\mathbb R}\).
\end{cor}

\begin{proof}
The projective limit of \(\eta(\tau)^{-1}\cdot[R_Z]\) as \(\tau\to0\) is represented by the sum of the terms of \(R_Z\) of maximal \(\eta\)-weight. 
By Theorem~\ref{thm:limit-cycle-decomposition} and multiplicativity of Chow forms, this sum is proportional to \(\prod_\sigma R_{Z_\sigma}\), with all Chow forms taken in \(\mathbb P(E^*)\). 
The weight polytope of the maximal-weight part is the exposed face, and the weight polytope of the product is the Minkowski sum of the weight polytopes of its factors.
\end{proof}

On \(X_\sigma\), the degree-one coordinates that survive are precisely those indexed by
\[
J_\sigma=\{\lambda\in J\cap\sigma:h_\lambda=g(\lambda)\}.
\]
Thus \(\pi_\sigma\) is defined by the linear system \(E_\sigma=\bigoplus_{\lambda\in J_\sigma}W_{1,\lambda}\).
Indeed, a summand first survives in \(\mathcal R_h/(\tau)\) exactly when \(h_\lambda=g(\lambda)\), and then survives in \(A_\sigma\) exactly when \(\lambda\in\sigma\). 

By Theorem~\ref{thm:limit-cycle-decomposition}, \(X_\sigma\) has moment polytope \(\sigma\) and open orbit \(G/H_{I_\sigma}\).
Its weight lattice is \(M\), its Duistermaat--Heckman polynomial is \(P_{DH}|_\sigma\), and its tail cone is \(\operatorname{cone}(I_\sigma)\). 

In the particular case when \(Y_\sigma=\mathbb P(E_\sigma^*)\), then \(Z_\sigma=d_\sigma[\mathbb P(E_\sigma^*)]\), where \(d_\sigma=\deg\pi_\sigma\). 
Its Chow form is the \(d_\sigma\)-th power of the determinant of the restrictions to \(E_\sigma^*\) of the \(n+1\) linear forms. 
Thus its Chow polytope is the singleton
\(
\left\{d_\sigma\sum_{\lambda\in J_\sigma}
\dim W_{1,\lambda}\,\varepsilon_\lambda\right\}
\).

If \(Y_\sigma\) has positive codimension in \(\mathbb P(E_\sigma^*)\), Theorem~\ref{thm:spherical-secondary-Chow} gives 
\[
\operatorname{Chow}_{\mathrm{GL}}^G(Z_\sigma)
=
(n+1)!\,\iota\left(
\int_\sigma K_{\sigma,x}P_{DH}(x)\,d\mu_M(x)
\right),
\]
where \(K_{\sigma,x}\) is the decomposition polytope from Proposition~\ref{def:decomposition-polytope}, formed with \(J_\sigma\) and \(\operatorname{cone}(I_\sigma)\).
We regard its points as vectors in \(\mathbb R^J\) by setting the coordinates outside \(J_\sigma\) to zero. 

Suppose now that \(X\) is horospherical. 
Since its tail cone is zero, \(\Delta=\operatorname{conv}(J)\).

\begin{defn}
For a triangulation \(\mathcal T\) of \(\Delta\) with vertices in \(J\), define its \emph{horospherical GKZ vector} \(\Phi_J^G(\mathcal T)\in\mathbb R^J\) by
\[
\Phi_J^G(\mathcal T)_\lambda
=
(n+1)!\!
\sum_{\substack{\sigma\text{ maximal in }\mathcal T\\
                 \lambda\in\operatorname{Vert}(\sigma)}}
\int_\sigma \beta_{\sigma,\lambda}(x)P_{DH}(x)\,d\mu_M(x),
\]
where \(\beta_{\sigma,\lambda}\) denotes the barycentric
coordinate on \(\sigma\) associated with the vertex \(\lambda\).
\end{defn}

If the heights \(h\) induce \(\mathcal T\) and the active weights on each maximal simplex are exactly its vertices, then Corollary~\ref{cor:Chow-face-decomposition} and the local barycentric formula give
\[
\operatorname{Face}_\eta
\bigl(\operatorname{Chow}_{\mathrm{GL}}^G(Z)\bigr)
=
\bigl\{\iota(\Phi_J^G(\mathcal T))\bigr\}.
\]

For arbitrary heights \(h\), the exposed face \(\operatorname{Face}_\eta(\operatorname{Chow}_{\mathrm{GL}}^G(Z))\) is the convex hull of the points \(\iota(\Phi_J^G(\mathcal T))\) for regular triangulations \(\mathcal T\) such that each maximal simplex is contained in a maximal cell \(\sigma\) and has all its vertices in \(J_\sigma\).

Indeed, the Chow polytope is the convex hull of the points \(\iota(\Phi_J^G(\mathcal T))\) for all regular triangulations: each such point is exposed by heights inducing the triangulation, and every vertex can be exposed by generic heights.
For any such \(\mathcal T\), let \(h_{\mathcal T}\) be the piecewise affine interpolation of \(h\) on its simplices. 
Then
\[
\bigl\langle\iota(\Phi_J^G(\mathcal T)),\eta\bigr\rangle
=(n+1)!\int_\Delta h_{\mathcal T}(x)P_{DH}(x)\,d\mu_M(x).
\]
Since \(h_{\mathcal T}\leq g_h\) and \(P_{DH}\) is positive on \(\operatorname{relint}(\Delta)\), this pairing attains the support value exactly when \(h_{\mathcal T}=g_h\). 
On each maximal simplex this means that \(g_h\) is affine and that its vertices satisfy \(h_\lambda=g_h(\lambda)\), giving the stated condition.
In particular, distinct regular triangulations give distinct vertices. 
Indeed, choose heights inducing \(\mathcal T\) for which the active weights on each maximal simplex are exactly its vertices. 
If another regular triangulation has the same horospherical GKZ vector, it must satisfy the face criterion for these heights. 
Each of its maximal simplices must therefore coincide with a maximal simplex of \(\mathcal T\), so the two triangulations coincide.

\begin{rem}
In the toric case, \(n=\dim\Delta\) and \(P_{DH}=1\).
The average of a barycentric coordinate on an \(n\)-simplex is \(1/(n+1)\), so
\[
\Phi_J^G(\mathcal T)_\lambda
=
\sum_{\substack{\sigma\text{ maximal in }\mathcal T\\
                 \lambda\in\operatorname{Vert}(\sigma)}}
n!\,\mu_M(\sigma).
\]
Thus \(\Phi_J^G(\mathcal T)\) is the classical GKZ vector~\cite[Chapter~7, Eq.~(1.4)]{GKZ}, with simplex volumes normalized by the lattice \(M\).
\end{rem}

\subsection{Examples}

\begin{exa}[A horospherical projective bundle]
\label{exa:horospherical-quadric-bundle}
Let \(Q=Q^3\) be the smooth quadric threefold, and consider \(X=\mathbb P_Q(\mathcal O\oplus\mathcal O(1))\) equipped with its anticanonical polarization \(L=-K_X\).
This is the rank-one horospherical Fano fourfold 4-1-50
of~\cite{Delcroix_Montagard}, under
\(G=\mathrm{Sp}_4\times\mathbb G_m\).

Let \(\omega_2\) denote the fundamental weight corresponding to \(Q\), and \(\chi\) the character of the additional torus.
Writing \(\beta=\omega_2+\chi\), we have \(M=\mathbb Z\beta\), \(\Delta=2\omega_2+[0,2]\beta\), and \(J=\{\lambda_0,\lambda_1,\lambda_2\}\) where \(\lambda_j=2\omega_2+j\beta\). 
Writing \(x=2\omega_2+s\beta\), the Duistermaat--Heckman measure is
\(
P_{\mathrm{DH}}\,d\mu_M
=\frac{(2+s)^3}{3}\,ds\). 

There are exactly two regular triangulations of \(\Delta\) with vertices in \(J\):
\[
\mathcal T_1=\{[0,1],[1,2]\},
\qquad
\mathcal T_2=\{[0,2]\},
\]
where we use the coordinate \(s\) to describe the simplices.

The horospherical GKZ vectors are
\[
\begin{aligned}
\Phi_J^G(\mathcal T_1)
=120\bigg(&
\int_0^1(1-s)P_{\mathrm{DH}}(s)\,ds,\\
&\int_0^1sP_{\mathrm{DH}}(s)\,ds
+\int_1^2(2-s)P_{\mathrm{DH}}(s)\,ds,\\
&\int_1^2(s-1)P_{\mathrm{DH}}(s)\,ds
\bigg)
=(262,1140,998),
\end{aligned}
\]
and
\[
\Phi_J^G(\mathcal T_2)
=120\left(
\int_0^2(1-s/2)P_{\mathrm{DH}}(s)\,ds,\,
0,\,
\int_0^2(s/2)P_{\mathrm{DH}}(s)\,ds
\right)=(832,0,1568).
\]
Consequently,
\[
\operatorname{Chow}_{\mathrm{GL}}^G(X)
=\iota\left(
\operatorname{conv}\{
(262,1140,998),(832,0,1568)\}
\right).
\]

Let us describe the corresponding degenerations.
For heights \(h=(0,1,0)\), the spherical envelope is \(g_h(s)=\min(s,2-s)\).
Its maximal cells are \([0,1]\) and \([1,2]\).
The corresponding components are embedded projective bundles \(Y_0=\mathbb P_Q(\mathcal O_Q(2)\oplus\mathcal O_Q(3))\) and \(Y_1=\mathbb P_Q(\mathcal O_Q(3)\oplus\mathcal O_Q(4))\), of degrees \(130\) and \(350\). 
The limiting Chow cycle is \([Y_0]+[Y_1]\).

For heights \(h=(0,-1,0)\), the spherical envelope is identically zero. 
There is a single maximal cell, but the middle block is inactive.
The normalized component is \(X\), and its active linear system defines a finite morphism of degree two onto
\(
Y=\mathbb P_Q(\mathcal O_Q(2)\oplus\mathcal O_Q(4)).
\)
The limiting Chow cycle is therefore \(2[Y]\), with \(\deg Y=240\).

Thus the two vertices of the equivariant Chow polytope correspond respectively to a reducible limit cycle and to a non-reduced limit cycle arising from a finite morphism.
\end{exa}

\begin{exa}[Chow vertices and non-simplicial subdivisions]
\label{exa:blowup-conic-P2P1-degenerations}
We return to the anticanonically polarized spherical
threefold \((X,L)\) of Example~\ref{exa:blowup-conic-P2P1},
whose spherical secondary polytope was computed there.
We retain the notation of that example and write
\(J=\{\lambda_0,\ldots,\lambda_4\}\), with the weights
ordered as before, and \(\varepsilon_0,\ldots,\varepsilon_4\)
for the corresponding characters of \(\widetilde T_0\).

We denote by \(\nu_1,\ldots,\nu_{10}\) the ten vertices
of the equivariant Chow polytope.
The following table gives their coordinates in the basis
\((\varepsilon_0,\ldots,\varepsilon_4)\), together with
representative heights exposing each vertex.
A height triple \((a,b,c)\) stands for \((a,b,c,0,0)\).

\[
\begin{array}{c|r|rrrrr}
 & (h_0,h_1,h_2)&p_0&p_1&p_2&p_3&p_4\\ \hline
\nu_1 &(-1,-1,-1)&18&0&0&88&54\\
\nu_2 &(-2,-1,1)&18&0&34&88&20\\
\nu_3 &(-1,1,2)&18&36&34&52&20\\
\nu_4 &(0,1,-1)&18&52&0&36&54\\
\nu_5 &(1,2,1)&18&52&18&36&36\\
\nu_6 &(3,2,-1)&34&36&0&20&70\\
\nu_7 &(4,3,1)&34&36&18&20&52\\
\nu_8 &(0,0,2)&36&0&52&52&20\\
\nu_9 &(2,0,-1)&52&0&0&20&88\\
\nu_{10}&(3,1,1)&52&0&36&20&52
\end{array}
\]
The corresponding subdivisions \(\mathcal D_i\) of the
moment polytope are pictured below, using \((t,s)\) as
planar coordinates.
The first five subdivisions refine the decomposition
along \(t=1\), whereas the remaining five do not.
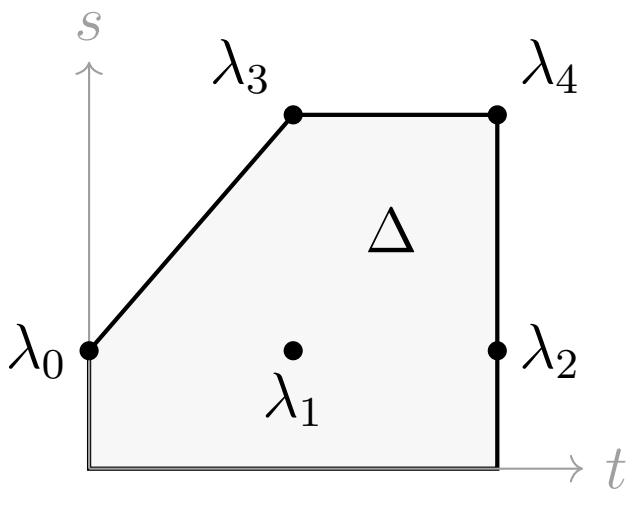
\begin{center}
\begin{tikzpicture}[x=1.35cm,y=.78cm]
\fill[gray!6] (0,0)--(0,1)--(1,3)--(2,3)--(2,0)--cycle;
\draw[thick] (0,0)--(0,1)--(1,3)--(2,3)--(2,0)--cycle;
\draw[->,gray!75] (0,0)--(2.42,0) node[right] {$t$};
\draw[->,gray!75] (0,0)--(0,3.45) node[above] {$s$};
\fill (0,1) circle (1.8pt) node[left] {$\lambda_0$};
\fill (1,1) circle (1.8pt) node[below] {$\lambda_1$};
\fill (2,1) circle (1.8pt) node[right] {$\lambda_2$};
\fill (1,3) circle (1.8pt) node[above left] {$\lambda_3$};
\fill (2,3) circle (1.8pt) node[above right] {$\lambda_4$};
\node at (1.48,2.02) {$\Delta$};
\end{tikzpicture}

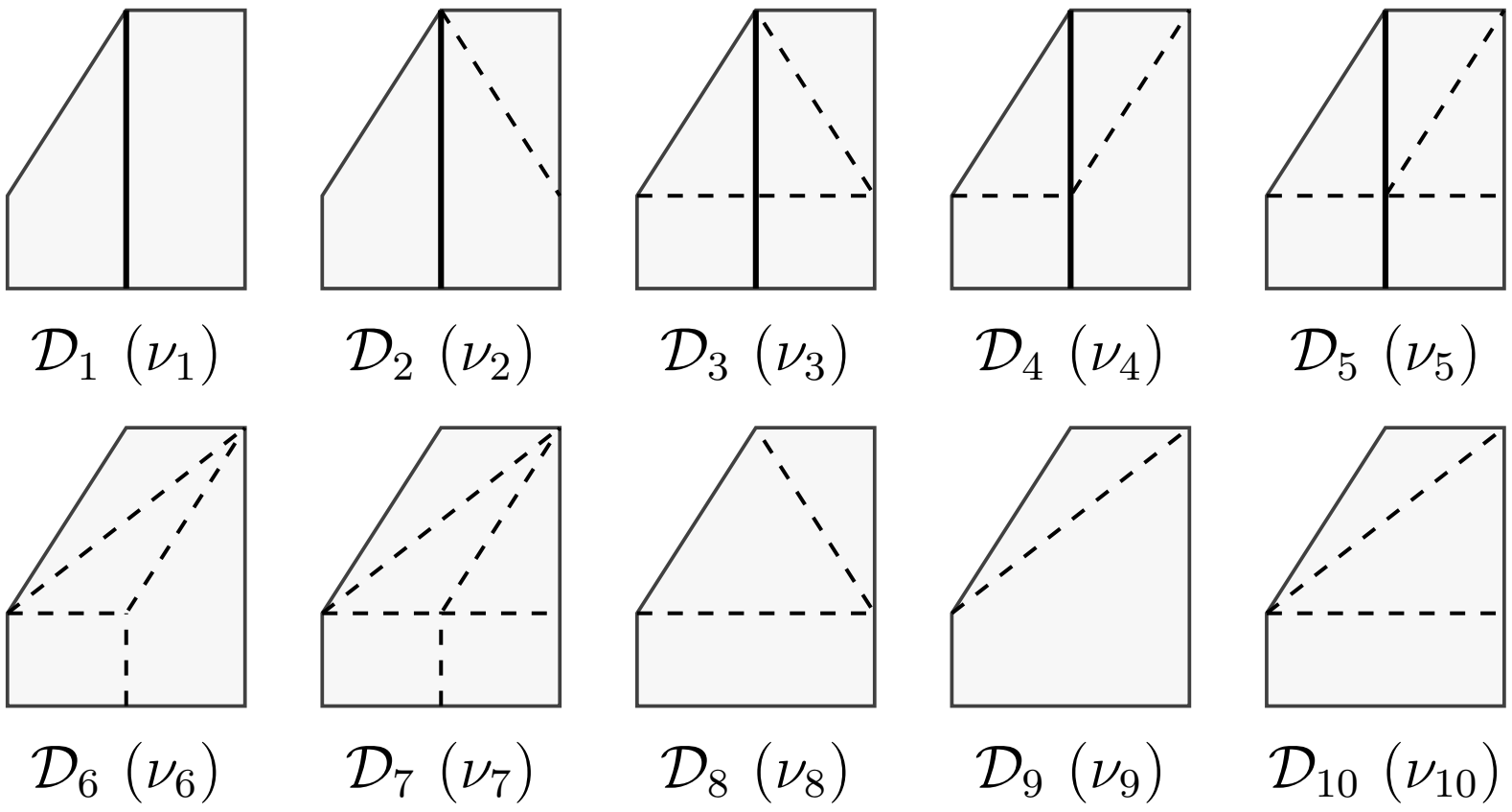
\medskip
\begin{tikzpicture}[x=0.82cm,y=0.64cm,
  wall/.style={draw=black, dashed, line width=.75pt},
  commonwall/.style={draw=black,line width=1.15pt}]
\newcommand{\SPSVsubdiv}[4]{%
  \begin{scope}[shift={(#1,#2)}]
    \fill[gray!6] (0,0)--(0,1)--(1,3)--(2,3)--(2,0)--cycle;
    \draw[black!75,line width=.7pt]
      (0,0)--(0,1)--(1,3)--(2,3)--(2,0)--cycle;
    #4
    \node[below] at (1,-.16) {$\mathcal D_{#3}\ (\nu_{#3})$};
  \end{scope}}
\SPSVsubdiv{0}{4.5}{1}{
  \draw[commonwall] (1,0)--(1,3);}
\SPSVsubdiv{2.65}{4.5}{2}{
  \draw[commonwall] (1,0)--(1,3);
  \draw[wall] (1,3)--(2,1);}
\SPSVsubdiv{5.3}{4.5}{3}{
  \draw[commonwall] (1,0)--(1,3);
  \draw[wall] (0,1)--(1,1)--(2,1)--(1,3);}
\SPSVsubdiv{7.95}{4.5}{4}{
  \draw[commonwall] (1,0)--(1,3);
  \draw[wall] (0,1)--(1,1)--(2,3);}
\SPSVsubdiv{10.6}{4.5}{5}{
  \draw[commonwall] (1,0)--(1,3);
  \draw[wall] (0,1)--(1,1)--(2,1);
  \draw[wall] (1,1)--(2,3);}
\SPSVsubdiv{0}{0}{6}{
  \draw[wall] (0,1)--(2,3);
  \draw[wall] (0,1)--(1,1)--(2,3);
  \draw[wall] (1,0)--(1,1);}
\SPSVsubdiv{2.65}{0}{7}{
  \draw[wall] (0,1)--(2,3);
  \draw[wall] (0,1)--(1,1)--(2,3);
  \draw[wall] (1,0)--(1,1)--(2,1);}
\SPSVsubdiv{5.3}{0}{8}{
  \draw[wall] (0,1)--(2,1)--(1,3);}
\SPSVsubdiv{7.95}{0}{9}{
  \draw[wall] (0,1)--(2,3);}
\SPSVsubdiv{10.6}{0}{10}{
  \draw[wall] (0,1)--(2,3);
  \draw[wall] (0,1)--(2,1);}
\end{tikzpicture}
\end{center}

In contrast with the toric case, these generic subdivisions
are not triangulations. In this example, the spherical
envelopes have vanishing \(s\)-slope on the quadrilateral
cells and strictly negative \(s\)-slope on the triangular
cells. By Theorem~\ref{thm:limit-cycle-decomposition},
the former correspond to non-horospherical components
with open orbit \(G/H\), and the latter to horospherical
components with open orbit \(G/H_{\varnothing}\).

Let us describe more precisely the degeneration corresponding
to \(\nu_9\), obtained for \(h=(2,0,-1,0,0)\).
By Lemma~\ref{lem:spherical-envelope-duality},  The spherical envelope is
\(
g_h(s,t)=\min(2-t,3-s),
\)
with two maximal cells
\[
\sigma_-=\Delta\cap\{s\leq1+t\},
\qquad
\sigma_+=\Delta\cap\{s\geq1+t\}.
\]
Their active weight sets and surviving spherical roots are
\[
\begin{array}{c|c|c}
 & J_\sigma & I_\sigma\\ \hline
\sigma_- & \{\lambda_0,\lambda_4\} & \{2\alpha\}\\
\sigma_+ & \{\lambda_0,\lambda_3,\lambda_4\} & \varnothing
\end{array}
\]
Thus the special fiber consists of two normal spherical components, one with the same open orbit as \(X\) and the other horospherical. 
Their degrees with respect to the pullback polarization are respectively \(26\) and \(14\). 
The Chow polytopes of the cycles \(Z_{\sigma_\pm}=(\pi_{\sigma_\pm})_*[X_{\sigma_\pm}]\) are
\[
\begin{aligned}
\operatorname{Chow}_{\mathrm{GL}}^G(Z_{\sigma_-})
&=\{36\varepsilon_0+68\varepsilon_4\},\\
\operatorname{Chow}_{\mathrm{GL}}^G(Z_{\sigma_+})
&=\{16\varepsilon_0+20\varepsilon_3+20\varepsilon_4\}.
\end{aligned}
\]
Their Minkowski sum recovers \(\{\nu_9\}\), as expected from Corollary~\ref{cor:Chow-face-decomposition}.

The chamber complex of the projection introduced in Remark~\ref{rem:BS},  as pictured in Example~\ref{exa:blowup-conic-P2P1},  is the common refinement of the ten subdivisions \(\mathcal D_i\). 
Note that the chamber complex has vertices which do not occur in any of the subdivisions \(\mathcal D_i\), for instance the point with coordinates \((t,s)=(4/3,7/3)\).
\end{exa}

\begin{exa}[Identical decompositions and distinct Chow vertices]
\label{exa:quadric-same-decomposition}
Consider the rank-two quadric of Example~\ref{exa:two-block-quadric}, with \(p=\dim U\), \(q=\dim W\), and \(X=\{q_U(u)+q_W(w)=0\} \subset\mathbb P(U\oplus W)\).
Recall that \(\Delta=\operatorname{conv}(0,\lambda_U,\lambda_W)\), \(J=\{\lambda_U,\lambda_W\}\), and that the spherical roots are \(2\lambda_U,2\lambda_W\).

For heights \(h=(h_U,h_W)\), the spherical envelope at \(x=s\lambda_U+t\lambda_W\) is
\[
g_h(x)=
\begin{cases}
(1-t)h_U+th_W,&h_U>h_W,\\
sh_U+(1-s)h_W,&h_U<h_W.
\end{cases}
\]
In both cases, the envelope is affine on all of \(\Delta\), and both degree-one weights are active.
Thus the decomposition consists of the single cell \(\Delta\), independently of the choice of generic heights.

Nevertheless, the corresponding degenerations differ. 
If \(h_U>h_W\), the slope annihilates \(2\lambda_U\) but not \(2\lambda_W\). 
The limiting variety is the singular quadric
\(
X_+=\{q_W(w)=0\}\subset\mathbb P(U\oplus W),
\)
with singular locus \(\mathbb P(U)\). 
Its open orbit is the satellite associated with \(I=\{2\lambda_U\}\).
Similarly, if \(h_U<h_W\), the limit is the singular quadric 
\(
X_-=\{q_U(u)=0\}\subset\mathbb P(U\oplus W),
\)
with singular locus \(\mathbb P(W)\) and surviving spherical root \(2\lambda_W\).

Thus the decomposition of the moment polytope, even together with its active weights, does not determine the exposed Chow face: the slopes of the spherical envelope, and hence the surviving spherical roots, must also be taken into account.
\end{exa}

\begin{exa}[A positive-dimensional Chow face with a single maximal cell]
\label{exa:GL-quadric-boundary}
This example shows that a degeneration with a single maximal cell may nevertheless expose a positive-dimensional local Chow polytope, even when all active weights are vertices of a simplex. 

Let \(m\geq2\), and consider the quadric 
\[
X=\{(u,\ell,z):\ell(u)+z^2=0\}
\subset\mathbb P\bigl(\mathbb C^m\oplus
(\mathbb C^m)^*\oplus\mathbb C\bigr),\]
with \(L=\mathcal O_X(1)\).
The group \(G=\mathrm{GL}_m\) acts by \(g\cdot(u,\ell,z)=(gu,\ell g^{-1},z)\).
Let \(a,b\) denote the highest-weights of the two nontrivial blocks of \(H^0(X,L)\).
Then \(M=\mathbb Za\oplus\mathbb Zb\), \(\Delta=\operatorname{conv}(0,a,b)\), and \(\Sigma=\{a+b\}\).

Choose block heights \(h_a=1\), \(h_b=-1\), \(h_0=0\) on the corresponding dual blocks.
The spherical envelope is \(g_h(sa+tb)=s-t\).
Its slope annihilates \(a+b\), so it lies in \(\operatorname{Lin}(\mathcal V)\). 
There is one maximal cell \(\sigma=\Delta\), all three blocks are active, and \(I_\sigma=\Sigma\). 
The degeneration is a product and its component map is the identity.

Nevertheless, for \(x=sa+tb\), the local decomposition polytope is the interval
\[
K_{\sigma,x}
=
\left\{(s+q,t+q,1-s-t-2q):
0\leq q\leq\frac{1-s-t}{2}\right\},
\]
in block order \((a,b,0)\).

Here \(n=2m-1\), \(d\mu_M=ds\,dt\), and
\[
P_{\mathrm{DH}}(sa+tb)
=
\frac{(s+t)s^{m-2}t^{m-2}}{(m-1)!(m-2)!}.
\]
The endpoints of \(K_{\sigma,x}\) are \((s,t,1-s-t)\) and \(\left(\frac{1+s-t}{2},\frac{1-s+t}{2},0\right)\).
Their integrals against \((2m)!\,P_{\mathrm{DH}}\,ds\,dt\) are respectively \((2m-1,2m-1,2)\) and \((2m,2m,0)\).
Since the intervals have the common direction \((1,1,-2)\), Corollary~\ref{cor:Chow-face-decomposition} gives, in full block coordinates,
\[
\operatorname{Chow}_{\mathrm{GL}}^G(X)
=
\operatorname{conv}
\{(2m-1,2m-1,2),(2m,2m,0)\}.
\]
\end{exa}

\section*{Acknowledgements}

This article is based on joint work by the authors over a period of more than three years. 
Many drafts were produced with different notation and different levels of generality. 
The authors made a significant use, in the Fall 2026, of ChatGPT to assist in putting together the various drafts, revising the exposition, reviewing and simplifying mathematical arguments. 
The authors reviewed the resulting material and take responsibility for the final manuscript. 

The first author thanks Sébastien Boucksom for discussions relating to this work. 

The first and second author were partially supported by ANR-21-CE40-0011 JCJC project
MARGE. 
The third author was partially supported by JSPS KAKENHI Grants JP24KK0252.

\bibliographystyle{alpha}
\bibliography{SPSV}

\end{document}